\documentclass[11pt,a4paper,english]{article}
\usepackage{import}
\usepackage{AHWpreamble}
\usepackage{caption}

\usepackage{natbib,bbm}
\usepackage[normalem]{ulem}
\usepackage{stackengine}
\newcommand\stackrqarrow[2]{%
    \mathrel{\stackunder[2pt]{\stackon[4pt]{$\approx$}{$\scriptscriptstyle#1$}}{%
            $\scriptscriptstyle#2$}}}

\renewcommand{\path}{\gamma}
\newcommand{\qsim}{\overset{\path}{\sim}}
\newcommand{\pathi}{q}
\newcommand{\hatgeod}{{\hat\path^*}}
\newcommand{\geod}{\path^*}
 \newcommand{\geodi}{\pathi^*}
\newcommand{\cpath}{\Gamma}
\newcommand{\cgeod}{{\Gamma^*}}

\newcommand{\cgeodi}[1]{{\Gamma^*_{#1}}}
\newcommand{\Rm}{\specialm g}

\newcommand{\length}{\mathtt{L}}
\newcommand{\cost}{\specialm{c}}

\title{Fermat costs and path-dependent first-passage percolation}
\date{}

\author{Sandro Gallo \\ Departamento de Estatística, Universidade Federal de São Carlos
        \and Matthieu Jonckheere \\   LAAS CNRS, Université de Toulouse
        \and Alejandro H. Wences \\ IIMAS, Universidad Nacional Autónoma de México}
\makeindex
\begin{document}
\maketitle
%\tableofcontents
\begin{abstract}
    We introduce a path-dependent and space-inhomogeneous first-passage percolation (FPP) model in the continuum. Inspired by Euclidean FPP, it is defined through a sample of points drawn from a common density function on a Riemannian manifold. In Euclidean FPP, the cost of each step in a discrete path is given by its length raised to the power $\alpha>1$. In our model, the power $\alpha$ increases as the path progresses, governed by a control mechanism. Consequently, the resulting percolation cost is not a distance. However, we prove a shape theorem showing that, as the number of points goes to infinity, the percolation costs converge to deterministic space-time integrals along curves in the manifold. This limit defines a quasi-distance with the property that traversing a curve from a high-density region towards a low-density one is generally cheaper than in the reverse direction; thus highlighting its potential use in machine learning applications.
    ~\\
   
\end{abstract}

\section{Introduction}

\paragraph{A space time inhomogeneous first passage percolation model} First-passage percolation (FPP), introduced by \cite{Hammersley}, provides a fundamental framework for studying random metric structures and optimal paths through disordered media. In the continuum setting, the seminal paper of \cite{HowardNewman1997} introduced the following percolation model:
\begin{equation}
\label{eq:defHN}
\forall x,y\in\Mani, \quad L_{n,\alpha}(x,y)\coloneqq \inf_{\substack{\path\colon x\overset{\path}{\sim}y\\ \gamma\subset \Xi\cup\{x,y\}}}\lbr \sum_{i=1}^{K-1} |\pathi_{i+1}-\pathi_i|^\alpha \rbr,
\end{equation}
where $\Xi$ is a unit-intensity Poisson process on $\mathbb R^d$, $\alpha>1$ and the infimum is taken over all possible discrete paths $\gamma$ that connect $x$ and $y$  through points in $\Xi$.
They established a shape theorem: the percolation distance
converges as $\abs{x-y}\to \infty$
to a deterministic multiple $\mu(d,\alpha)$ of Euclidean distance, the
constant depending on the environment only through $d$ and $\alpha$. The
superadditivity underlying \eqref{eq:defHN} rests on two features: the cost of a
step is a function of its length alone, and the environment is statistically
invariant under translation.\par 
The second feature has been relaxed in the spatial direction by \cite{HwangDamelinHero2016} and  \cite{GroismanJonckheereSapienza2022}. They replaced
the homogeneous Poisson process by an i.i.d.\ sample $Q_n$ of $n$ points
drawn from a density $f$ on a compact $d$-dimensional Riemannian manifold
$\Mani$ embedded in $\mathbb R^D$; and studied the (microscopic) sample Fermat distance:
\begin{equation}
\label{eq:defFDclassicalIntro}
%\label{sample.fermat}
\forall x,y\in\Mani, \quad D^\dist_{n,\alpha}(x,y)\coloneqq \inf_{\substack{\path\colon x\overset{\path}{\sim}y\\ \gamma\subset Q_n\cup\{x,y\}}}\lbr \sum_{i=1}^{K-1} \dist(\pathi_{i+1},\pathi_i)^\alpha \rbr,
\end{equation}
where $\dist$ is either the Riemannian distance $\dist_\Mani$, or the Euclidean distance $\dist_E$ in the ambient space $\R^D$.
After scaling by $n^{(\alpha-1)/d}$ this converges to $\mu(\alpha,d)D_{\infty,\alpha}(x,y)$ where $D_{\infty,\alpha}(x,y)$ is the (macroscopic) Fermat distance:
\begin{equation}
 \label{def:FD}
\forall x,y\in\Mani,\quad D_{\infty,\alpha}(x,y)\coloneqq \inf_{\cpath\colon x\overset{\cpath}{\sim}y} \int_0^{\length\Gamma} f^{-\frac{\alpha-1}{d}}(\cpath(t)) dt.
\end{equation}
The infimum is taken over piece-wise smooth curves $\cpath$, with length $\length\cpath>0$, that are parameterized with speed 1.  The medium is now
inhomogeneous, but the inhomogeneity is purely \emph{spatial},  in particular the limit remains a length functional. More precisely, \eqref{def:FD} defines a geodesic distance in a metric which is conformal to $\Rm$, the Riemannian metric on $\Mani$. This fact is also what limits the use of
\eqref{eq:defFDclassicalIntro} as a tool for learning a geometry from data, we will come back to this point later on.\par
In this work we study first-passage percolation in a medium that is
inhomogeneous in space \emph{and} in time: the power parameter $\alpha$ in \eqref{eq:defFDclassicalIntro} now increases as
the path progresses (see Section \ref{sec:MainDefinitions} for formal definitions). The macroscopic object, the Fermat cost, is then no longer
a length functional but a space-time integral along the path:
\begin{equation}\label{eq:macroFermatCost}
C_{\infty}(x,y)\equiv C_{\infty,\alpha,A}(x,y)\coloneqq \inf_{\substack{\cpath\colon x\overset{\cpath}{\sim}y}} \int_0^{\length\cpath} e^{-\frac{A(t)}{d}} f^{-\frac{\alpha-1}{d}}(\cpath(t)) dt.
\end{equation}
The function $A\colon[0,\infty)\to[0,\infty)$, assumed bounded, non-decreasing, and uniformly continuous, is a parameter of the model.
It will prescribe the evolution of the power parameter $\alpha$ along a path in the microscopic percolation mechanism, modeling how the local cost evolves as time progresses, as measured by the distance already traveled. 

\paragraph{Learning applications and the asymmetric Fermat cost}

Meaningful distances are rarely available as direct inputs for machine learning tasks like clustering or recommendation, and standard Euclidean or feature-space distances in high-dimensional space degrade due to distance concentration. The ISOMAP distance \cite{Tenenbaum2000} has been one of the most popular methods for nonlinear dimensionality reduction, but it doesn't take into account the density of sample points. As an alternative \cite{SapienzaGroismanJonckheere2018,GroismanJonckheereSapienza2022} suggested to use the distance functions in \eqref{def:FD} that bring together points which can be joined
by paths that go through high-density regions of $\Mani$. It appears also that sample Fermat distances \eqref{eq:defFDclassicalIntro} perform well as estimators of \eqref{def:FD} in practice \cite{SapienzaGroismanJonckheere2018,GroismanJonckheereSapienza2022,McKenzie2019, HwangDamelinHero2016,FernandezBorghiniMindlinGroisman2023}.
\par 
What \eqref{def:FD} cannot do is exactly what the time-inhomogeneity in \eqref{eq:macroFermatCost} provides.
The distance \eqref{def:FD} is symmetric, so it can record how far apart two points are but
not how the landscape between them is arranged. In contrast \eqref{eq:macroFermatCost} is asymmetric: a path is in general cheaper to cross from a high-density region towards a low-density one, than in the reverse direction. In practice, the sign of $C_\infty(x,y) - C_\infty(y,x)$
might provide a useful density ranking of the two end points $x$
and $y$.
We expect this additional features to be exploitable in a range of learning tasks, including clustering, outlier detection, and recommendation, opening new directions in machine learning that symmetric distance-based methods cannot access.
Importantly, this is achieved without having to estimate the density $f$ directly, thus avoiding the curse of dimensionality \cite{Vershynin2018, hastie2009elements}.

\paragraph{Main mathematical challenges to reach \eqref{eq:macroFermatCost}}

We highlight three main differences that are determinant to our proofs in the study of \eqref{eq:macroFermatCost} when compared to  $D_{\infty,\alpha}$.\par
    \emph{Loss of symmetry:} The weight $e^{-A(t)/d}$ decreases along the
path and $C_\infty(x,y)\neq C_\infty(y,x)$ in general. The limit is a cost function, not a
distance. Convergence nevertheless takes place in a Gromov-Hausdorff-type
topology, provided one works with \emph{cost spaces} (compact networks in the
sense of \cite{ChowdhuryFacundo2023}). This is Theorem \ref{th:GHconvergence2}.\par
    \emph{Loss of subpath optimality:} A subsegment of an optimal curve for
\eqref{eq:macroFermatCost} need not be optimal for its own endpoints. For instance
a detour of a path can pay for itself by buying a larger discount factor $e^{-\frac{A(t)}{d}}$
to traverse expensive ground in the future (small values of $f$). This is an important
technical obstacle of the paper. Section \ref{subsec:locpropcinfty} supplies the replacement: the
infimum in \eqref{eq:macroFermatCost} may be restricted to a class $\mathcal
F_{\kappa^*}$ of curves of uniformly bounded curvature and length, thus avoiding geodesic paths with deep contortions. \par
    \emph{The cost of advancing:} The variable $t$ in
\eqref{eq:macroFermatCost} is arc-length along the macroscopic curve, and the
microscopic model has to manufacture it out of the discrete approximation. The obvious
choice, advancing the exponent $\alpha$ with the accumulated length traversed by a discrete path, fails: a path can
accumulate length by shuttling between two very close sample points,
so the clock saturates reaching a maximum value of $\alpha$ at a negligible cost, and the whole space-time structure
collapses. To avoid the latter we introduce a control mechanism into the percolation dynamics. In this model the power parameter $\alpha$ increases along the path whenever a marked point in the path is reached. The control mechanism consists of choosing the appropriate marks within a path, under the restriction that marks must be at a distance larger than a mesoscopic distance threshold $\delta_n$ (see Definition \ref{def:FermatCostGeneral}). This ensures that a true displacement of $\delta_n$ is made at each change point of $\alpha$, so that a cost of at least $\approx\mu(\alpha,d)\maxv f^{-\frac{\alpha-1}{d}}$ is incurred. Driving the clock to saturation now costs an amount bounded away from zero, so stalling is still possible but priced.

\paragraph{Organization of the paper}
Section~\ref{sec:main} introduces the microscopic FPP models and states the main results: the Shape Theorem for microscopic Fermat costs  (Theorem~\ref{th:FermatCostConvergence}) and the convergence of cost spaces (Theorem~\ref{th:GHconvergence2}). 
%and the Shape Theorem for Euclidean path-dependent FPP (Theorem~\ref{th:DynamicalEFPP}). 
In section \ref{subsec:proofOverview_heuristics} we give heuristics and overall explanations for the proof, as well as a description of its layout into the remaining sections of the manuscript.

\section{A path-dependent FPP model: Main results}
\label{sec:main}

\subsection{Definition of the microscopic model}\label{sec:MainDefinitions}
Consider a Riemannian manifold $\Mani$ of dimension $d\geq1$, with geodesic metric $\dist_\Mani$, and embedded in the Euclidean space $\R^D$, $D\geq d$. As before, when we only use $\dist$ as index, we mean that the corresponding quantity is defined using either 
the Euclidean $\dist_E$ of the Riemannian $\dist_\Mani$ distance. 
Let 
$A\colon \Rp\to[0,\infty)$ be a bounded non-decreasing function and $0<\delta_n\leq 1$, a sequence of positive  numbers satisfying
\begin{align}
    \delta_n&\toas{n\to\infty}0,
%    \nonumber\\
%    \frac{\log(n)}{n\delta_n^d} &\toas{n\to\infty}0.
\end{align}
(see Assumption \ref{cond:DeltanBound} below for further assumptions on $\delta_n$).\par
Fixing  $\alpha>1$, we study the following extension of \eqref{eq:defFDclassicalIntro} into a path-dependent setting with control. For any discrete path $\path=(\pathi_1,\cdots,\pathi_{\card \gamma})\subset \Mani \subset \R^D$ 
we define the set $M_{\dist,\delta_n}(\path)$ of $\dist$-admissible marks $(\iota(k))_{k=1}^{\ell}\equiv(\iota_\path(k))_{k=1}^{\ell}$  of parameter $\delta_n$ by
\begin{align*}
M_{\dist,\delta_n}(\path) \coloneqq \{(\iota(k))_{k=1}^{\ell}\colon &1\leq \ell \leq \card\path,\quad \iota(1)=1, 
\\&\forall 1\leq k \leq \ell-1, \quad \dist(\pathi_{\iota(k)},\pathi_{\iota(k+1)})\geq\delta_n \}.
\end{align*}

The key idea of the model is that the cost of each step along the path is updated by increasing the power parameter $\alpha$ in \eqref{eq:defFDclassicalIntro}, as the path traverses the different marks. 
Specifically, within the path segment between two consecutive marks, say $(\pathi_{\iota(k)},\cdots,\pathi_{\iota(k+1)})$, the power parameter is set to 
$$
	\alpha_k \equiv \alpha_{n,k} \coloneqq \alpha + \frac{A(\delta_n k)}{\log(n)}
$$
% $\alpha+ \frac{A(\delta_n k)}{\log(n)}$ 
(where $\alpha\ge1$) so the cost of a step $(\pathi_i,\pathi_{i+1})\subset \path$, \emph{within this segment}, is 
$$
\dist(\pathi_{i+1},\pathi_i)^{\alpha+ \frac{A(\delta_n k)}{\log(n)}}.
$$
Letting
$$
	I^\iota_\path(i) \equiv I^{\iota}_{n,\path}(i) \coloneqq k \iff \iota_\path(k)\leq i <\iota_{\path}(k+1)
	$$
(with the convention that $\iota(\ell+1)=+\infty$) we can write down the cost of a generic step $(\pathi_i,\pathi_{i+1})\subset \path$ as
    $$
\dist(\pathi_{i+1},\pathi_i)^{\alpha_{I^{\iota}_{\path}(i)}}.
$$

We refer the reader to the \say{Scaling Heuristics} paragraph in section \ref{sec:MainResults}, as well as to section \ref{subsec:proofOverview_heuristics}, for a heuristic argument explaining the specific form of the $\alpha_{n,k}$.

%The above leads to the following formal definition (see also Figure \ref{fig:pathDependentFPP}).
\begin{definition}[Microscopic $\dist$-Fermat Cost]\label{def:FermatCostGeneral}
 The   \emph{microscopic $\dist$-Fermat Cost} between a pair of points $x,y\in\Mani$ is defined as
\begin{equation}\label{eq:defDFermatMicroGeneral}
 C^\dist_n(x,y)\coloneqq \inf_{\substack{\path\colon x\overset{\path}{\sim}y\\ \path\subset Q_n \cup \{x,y\} }}
 \inf_{\iota\in M_{\dist,\delta_n}(\path)}
 \sum_{i=1}^{\card \path -1}
\dist(\pathi_{i+1},\pathi_i)^{\alpha_{I^\iota_\path(i)}}.
\end{equation}
\end{definition}

$C_n^\dist$ will be called Euclidean (\emph{resp.} Geodesic) Fermat cost if $\dist=\dist_E$ (\emph{resp.} if $\dist=\dist_\Mani$). It is important to keep in mind that the indices $\iota_\path(k)$ used in the definition of $C_n^\dist$ also depend on the distance $\dist$ which is used. 
 
\begin{remark}
    Note that if $A\equiv0$ then $C^\dist_n(x,y)=D^\dist_{n,\alpha}(x,y)$.
\end{remark}
\begin{remark}
    Note that $C^\dist_n(x,y)\ne C^\dist_n(y,x)$, and also the triangle inequality is in general not satisfied since after concatenating two marked paths, the marks in the resulting path may not be admissible. For this reason we refer to $C^\dist_n$  as a cost and not as a distance nor a quasi-distance.\par
    The macroscopic Fermat cost \eqref{eq:macroFermatCost} does, however, satisfy the triangle inequality. This is an easy consequence of the fact that the function $A$ is non-decreasing, so that the cost of traversing the concatenation
    of two curves is smaller than the sum of the separate costs. Thus $C_\infty$ is a quasi-distance.
\end{remark}

\subsection{Assumptions}\label{sec:assumptions}
For the reader's convenience, we compile the main assumptions needed to state our results.

\begin{assumption}\label{cond:continuousMetricSpace}
    The Riemannian manifold $\Mani$ is smooth and compact. In particular 
    $$
        \diameter(\Mani):=\max_{x,y\in\Mani}\dist_\Mani(x,y)<\infty.
    $$
    Also, the manifold $(\Mani,\Rm)$ is embedded in the Euclidean space $\R^D$,  and we assume that the embedding has bounded geometry, i.e. that the norm of its second fundamental form $\specialm{II}$ is uniformly bounded on $\Mani$ by some constant $\const_{\specialm {II}} < \infty$. (Note that this is always satisfied if $\Mani$ is smoothly embedded in $\R^D$ as $\Mani$ is compact).
\end{assumption}
\begin{assumption}\label{cond:boundeddensity}
    The density $f$ of $\mu$ is smooth and has a bounded gradient. Also it satisfies
    $$0\le \minv f \coloneqq \inf_{x\in\Mani}f(x) \leq \maxv f\coloneqq \sup_{x\in\Mani}f(x)\le \infty.$$
\end{assumption}
\begin{assumption}\label{cond:functionA}
    The power parameter $\alpha$ satisfies $$\alpha>1.$$\par 
    Also, the function $A\colon [0,\infty)\to[0,\infty)$ is  uniformly continuous, non-decreasing, and bounded:
    $$
        0\leq \maxv A \coloneqq \sup_{t\in [0,\infty)}A(t) \le\infty.
    $$
\end{assumption}
The following assumption ensures that $\delta_n$ is of mesoscopic scale, allowing us to obtain Riemann–Stieltjes sum approximations of the 
curve integrals appearing in \eqref{eq:macroFermatCost}, while also ensuring that whenever $\dist_\Mani(q,q')\geq\delta_n$ we have $n^{\frac{\alpha-1}{d}}D_{n,\alpha}(q,q')\approx \mu(d,\alpha)D_{\infty,\alpha_{n,k}}(q,q')$ (see sections \ref{subsec:proofOverview_heuristics} and \ref{sec:macroscopicLimit}).  
\begin{assumption}\label{cond:DeltanBound}
    There exist $0<\eta<\eta'<1/\alpha$ and $C,C'>0$ such that the radii $\delta_n$ also satisfy 
    $$
     C n^{-\frac{1-\eta}{d}}\leq \delta_n \leq C' n^{-\frac{1-\eta'}{d}}.
    $$
\end{assumption}

As in \cite{HwangDamelinHero2016}, our convergence results for microscopic Fermat costs will hold uniformly for all pairs of points $x,y\in\Mani$ which are not too close to each other, that is, such that $\dist(x,y)\geq b_n$ where $b_n$ is allowed to go to zero, but not too fast:
\begin{assumption}\label{cond:bn}
    For $0<\eta'<1/\alpha$ as in Assumption \ref{cond:DeltanBound}, the minimum distances $(b_n)_{n\geq 1}$ between points $x,y\in\Mani$ satisfy
	$$
	\frac{n^{\frac{\alpha\eta'-1}{d}}}{b_n}\toas{n\to\infty} 0.
	$$
\end{assumption}

\subsection{Statements of the main results}\label{sec:MainResults}
Let  $Q_n=\{q_1,\ldots,q_n\}\subset \Mani$ be a sample of $n\geq 1$, i.i.d. points drawn from a common probability measure $\mu$ defined on $\Mani$ with density $f$.
Recall Definition \ref{def:FermatCostGeneral} of the Fermat cost $C^\dist_n(x,y)$, with $\dist$ being either $\dist_\Mani$ or $\dist_E$, between any pair of points in $x,y\in\Mani$. 

\paragraph{Equivalence between $C^{\dist_\Mani}_n$ and $C^{\dist_E}_n$} The functional $C^{\dist_E}_n$ is more natural from an applied perspective: in distance learning, the manifold distance $\dist_\Mani$ is typically unavailable, so $C^{\dist_E}_n$ serves as a practical estimator. A similar observation was made by \cite{GroismanJonckheereSapienza2022} in the context of the Fermat distance. On the other hand, when the primary interest is mathematical, the functional $C^{\dist_\Mani}_n$ is more natural and appears easier to manipulate.  Our first result relates these two perspectives by establishing that both functionals are asymptotically equivalent. As a consequence, any result proved for one automatically transfers to the other. This is true for the results of this paper and in particular also for those of \cite{HwangDamelinHero2016} and \cite{GroismanJonckheereSapienza2022}. 

\begin{proposition}\label{prop:equivalentModels}
	Fix $\eps>0$. There exist constants $\const_0>0,\const_1>0$ 
	such that for large enough $n$,
	\begin{equation}
		\P\left(\sup_{\substack{x,y\in\Mani\\\dist_\Mani(x,y)\geq b_n}}  \abs{\frac{C^{\dist_E}_{n}(x,y)}{C^{\dist_\Mani}_n(x,y)} -1} \ge \eps \right) \leq e^{-\const_0n^{\const_1}}. 
	\end{equation}
\end{proposition}

\noindent{\bf Notation alert.} Owing to the above proposition, we will state some results using $C^\dist_n$, meaning that either distance can be used. 

\paragraph{Shape Theorem for Fermat Costs}	
 The main result of our work is the following path-dependent Shape Theorem which holds under the assumptions stated in section \ref{sec:assumptions}.
\begin{theorem}\label{th:FermatCostConvergence}
     Fix $\eps>0$. There exist constants $\const_0>0,\const_1>0$ 
    such that for large enough $n$,
    \begin{equation}
            \P\left(\sup_{\substack{x,y\in\Mani\\\dist(x,y)\geq b_n}}  \abs{\frac{n^{\frac{\alpha-1}{d}}C^\dist_{n}(x,y)}{\mu(d,\alpha)C_\infty(x,y)} -1} \ge \eps \right) \leq e^{-\const_0n^{\const_1}},
    \end{equation}
    where, we recall, $C_{\infty}(x,y)$ was defined in \eqref{eq:macroFermatCost}. In particular, the following convergence holds almost surely for fixed $x,y\in\Mani$,
    $$n^{\frac{\alpha-1}{d}}C^\dist_{n}(x,y)\stackrel{n\to\infty}{\longrightarrow} \mu(d,\alpha)C_\infty(x,y).$$
    
\end{theorem}
The limit $\mu(d,\alpha)C_\infty(x,y)$ is the same,  regardless of which distance, $\dist_\Mani$ or $\dist_E$, is used in the definition. 

\paragraph{Scaling heuristics.}
Let $\Delta_n$ be the \emph{spacing}  of the sample $Q_n$ on $\Mani$, which is defined as 
\begin{equation}\label{eq:defSpacing}
\Delta_n\coloneqq \sup_{x\in\Mani}\dist_\Mani(x,Q_n).
\end{equation}
The spacing satisfies $\log(\Delta_n) \sim -\frac{1}{d}\log(n)$ and is the minimal scale at which $Q_n$ resolves the geometry of $S$; forcing steps in Fermat geodesics to be of order $n^{-\frac{1}{d}}$ is therefore natural (see Lemma \ref{le:DeltanBound}). In return, it forces the step cost at $\ell^\alpha \sim n^{-\alpha/d}$ and the number of steps at $K_n \sim n^{1/d}$, giving a total cost $C^\dist_n \sim n^{-(\alpha-1)/d}$. As a consequence, the scaling $n^{(\alpha-1)/d}$ is thus the unique polynomial normalisation making the rescaled cost $O(1)$.  Any other polynomial scaling would make this single-step contribution either vanish or blow up, forgetting the local density information. The exponent perturbation $A(k\delta_n)/\log n$ is added without disturbing this balance: since $\log\ell \approx -(\log n)/d$, one gets $\ell^{A(k\delta_n)/\log n} \to e^{-A(k\delta_n)/d}$, contributing only a bounded, path-dependent multiplicative weight per step. The two requirements — steps at the smallest admissible spacing, each step carrying local density information — therefore justify the scaling $n^{(\alpha-1)/d}$ and the $1/\log n$ form of the exponent perturbation.

%\sandro{I suggest to put the scaling heuristic here.}

\paragraph{Convergence of cost spaces}
The previous Theorem leads to an extension of the convergence result of the main Theorem in \cite{FernandezBorghiniMindlinGroisman2023} (see section 1.2 therein).
They prove the convergence, in the Gromov-Hausdorff topology, of the random metric space $(Q_n,n^{\frac{\alpha-1}{d}}D^{\dist_E}_{n,\alpha}(\cdot,\cdot))$ towards the deterministic
metric space $(\Mani,\mu(d,\alpha)D_{\infty,\alpha}(\cdot,\cdot))$. The latter has important implications for applications in topological data analysis such as the convergence of persistence homology diagrams, one of the main tools used in this area (see therein for more details). We extend the above convergence of metric spaces into our setting in which we have cost functions rather than distances. 

The Gromov-Hausdorff distance $\dist_{GH}$ between two metric spaces $(\Mani_1,\dist_1)$ and $(\Mani_2,\dist_2)$  (see e.g. Definition 3.4 in \cite{Gromov2007})  has an equivalent
formulation in terms of the so called correspondence sets $A$
% an their distortions 
(see e.g. section 7.3 \cite{BuragoBuragoIvanov2001}). Correspondence sets are subsets $A\subset\Mani_1\times\Mani_2$ whose projections on $\Mani_1$ and $\Mani_2$ satisfy
 $\pi_{\Mani_1}(A)=\Mani_1$
and $\pi_{\Mani_2}(A)=\Mani_2.$ 
% The distortion of a correspondence set $A$ is given by $\sup_{(x,y),(x',y')\in A}
% \abs{\dist_1(x,x') - \dist_2(y,y')}$. \sandro{we don't use the notion of distortion nor its notation...}
The Gromov-Hausdorff distance between $(\Mani_1,\dist_1)$ and $(\Mani_2,\dist_2)$ can then be written as:
$$
\dist_{GH}((\Mani_1,\dist_1),(\Mani_2,\dist_2)) = \frac{1}{2}\inf_{A\subset \Mani_1\times\Mani_2} \sup_{(x,y),(x',y')\in A}
\abs{\dist_1(x,x') - \dist_2(y,y')},
$$
where the infimum is taken over all correspondence sets $A$.\par
Define a \emph{cost space} $(\Mani,\cost)$ as a topological space $\Mani$ that is endowed with a positive continuous cost function $\cost\colon \Mani\times\Mani\to [0,\infty)$ between its elements, that satisfies $\cost(x,y)=0$ if and only if $x=y$. Examples of cost spaces include quasi-metric spaces for which $\cost$ satisfies the triangular inequality, and metric spaces for which $\cost$ is a distance function. In turn, cost spaces are an important example of a network, as defined in \cite{ChowdhuryFacundo2023}. Two cost spaces $(\Mani_1,\cost_1)$
and $(\Mani_2, \cost_2)$  are said to be weakly isomorphic if there exists a set $Z$ and two surjective maps $\phi_{\Mani_1}\colon Z\to\Mani_1$ 
and $\phi_{\Mani_2}\colon Z\to\Mani_2$ such that 
$\cost_1(\phi_{\Mani_1}(z),\phi_{\Mani_1}(z'))=\cost_2(\phi_{\Mani_2}(z),\phi_{\Mani_2}(z'))$ for all $z,z'\in Z$. Inspired by $\dist_{GH}$ above, the function $\dist_{\mathcal N}$ below defines a distance function on the class of compact cost spaces, and more generally on the class of compact networks, up to weak isomorphism \cite[Theorem 12]{ChowdhuryFacundo2023}:
$$
\dist_{\mathcal{N}}((\Mani_1,\cost_1),(\Mani_2,\cost_2)) \coloneqq \frac{1}{2}\inf_{A\subset \Mani_1\times\Mani_2} \sup_{(x,y),(x',y')\in A}
\abs{\cost_1(x,x') - \cost_2(y,y')},
$$
where again the infimum is taken over all correspondence sets $A$. 
We prove the following result.
\begin{theorem}\label{th:GHconvergence2}
    For any fixed $\eps>0$, there exist constants $\const_0>0,\const_1>0$ 
	such that for large enough $n$,
	$$
	\P\left(  \dist_{\mathcal{N}}((Q_n, C^\dist_n), (\Mani, \mu(d,\alpha)C_\infty)) > \eps \right) \leq e^{-\const_0n^{\const_1}}.
	$$
\end{theorem}
\par

\subsection{Proof overview and structure}
\label{subsec:proofOverview_heuristics}

\paragraph{Strategy of proof.}
The central result, Theorem~\ref{th:FermatCostConvergence}, asserts that
$$
  n^{\frac{\alpha-1}{d}}\,C_n^\dist(x,y)
  \;\approx\;
  \mu(d,\alpha)\,C_\infty(x,y),
$$
that is, the properly rescaled microscopic (random) cost converges to the macroscopic (deterministic) cost. In order to prove this, we introduce an intermediate \emph{mesoscopic} cost $\widehat{C}_{n,\infty}^\dist(x,y)$ (defined in Section~\ref{sec:mesoscopicApproximation}) and split the convergence into two steps:
$$
n^{\frac{\alpha-1}{d}}\,C_n^\dist(x,y)  \;\stackrqarrow{\text{Theorem \ref{theo:hatCnVshatCninfty}}}{\text{Section \ref{sec:mesoscopicApproximation}}}\;
\mu(d,\alpha)\,\widehat{C}_{n,\infty}^\dist(x,y)
  \;\stackrqarrow{\text{Theorem \ref{th:hatCnInfVSCInf}}}{\text{Section \ref{sec:macroscopicLimit}}}\;
  \mu(d,\alpha)\,C_\infty(x,y).
$$

\begin{itemize}
  \item \textbf{Step~1: Microscopic to mesoscopic} (Section~\ref{sec:mesoscopicApproximation}). We group the small-scale percolation steps into mesoscopic segments, each of length of order $\delta_n$. Within each segment, the power parameter is essentially constant, so the existing shape theorem for standard Fermat distances (Proposition~\ref{prop:alphanDnConv}) can be applied segment by segment. Concatenating these local approximations yields Theorem~\ref{theo:hatCnVshatCninfty}.
  \item \textbf{Step~2: Mesoscopic to macroscopic} (Section~\ref{sec:macroscopicLimit}). The mesoscopic cost $\widehat{C}_{n,\infty}^\dist(x,y)$ is a discrete sum of the form $\sum_k e^{-A(k\delta_n)/d}\,D_{\infty,\alpha}(u_k,u_{k+1})$, where the $u_k$ are mesoscopic waypoints. As $\delta_n\to 0$, this sum is a Riemann approximation to the space-time path integral defining $C_\infty(x,y)$. Making this precise is the content of Theorem~\ref{th:hatCnInfVSCInf}.
  \item 
  \textbf{Step~3: Concluding}
{(Section~\ref{sec:final_step}).} We combine the preceding steps to prove the main results: Theorems~\ref{th:FermatCostConvergence} and  \ref{th:GHconvergence2}.
\end{itemize}

\paragraph{Organisation of the rest of the manuscript.}
Section~\ref{sec:prelim} collects notation, definitions, and auxiliary results used throughout; the proofs of these auxiliary results are given in Appendix~\ref{sec:appendix}. Sections \ref{sec:mesoscopicApproximation}, \ref{sec:macroscopicLimit} and \ref{sec:final_step} are, as explained above, the three main steps toward the proof of the theorems. Proposition~\ref{prop:equivalentModels}, which establishes the asymptotic equivalence of the Euclidean and geodesic Fermat costs $C_n^{\dist_E}$ and $C_n^{\dist_\Mani}$, allows us to work with whichever distance is more convenient. In particular, in Sections \ref{sec:mesoscopicApproximation} and \ref{sec:final_step}, which involve microscopic costs, we use the geodesic distance and write simply $C_n = C_n^{\dist_\Mani}$.  
Proposition~\ref{prop:equivalentModels} itself is proved in in appendix \ref{sec:step0}. 

\section{Preliminaries and further notation} \label{sec:prelim}
This section contains several technical results that we consider auxiliary to prove our theorems. Their proofs are postponed to Appendix \ref{sec:appendix}.
Before we state them, let us introduce some further notation.
\subsection{Notation}
 We will write 
$b\pm c$ to denote a value $a\in\R$ (typically not explicit) that satisfies $\abs{a-b}\leq c.$ %\par 
Also we will write
$g(x)=\lo(h(x))$ as $x\to a$ whenever $\lim_{x\to a}g(x)/h(x)=0$. We will say that
$g(x,y)=\lo(h(x))$ uniformly on $y\in Y$ whenever  $\lim_{x\to a}g(x,y)/h(x)=0$ uniformly on $y\in Y$.
Similarly we will write $g(x)=\bO(h(x))$ as $x\to a$ whenever $\lim_{x\to a}g(x)/h(x)\leq C$ for some $C>0$,
and uniformly on $y\in Y$ whenever $\sup_{y}g(x,y)/h(x)\leq C$ for all $x$ close enough to $a$. For a function $h\colon X\to\R$ defined on a metric space $(X,\dist_X)$ let 
$\modulus_h\colon \Rp \to \Rp$ be its modulus of continuity. That is 
$$
\modulus_h(\delta) \coloneqq \sup_{\substack{x,x'\in X\\ \dist_X(x,x')\leq\delta}} \abs{h(x)-h(x')}.
$$
Note that $h$ is uniformly continuous if and only if $\modulus_h(\delta)\to0$ as $\delta\to0$.
\par
Recall that the microscopic Fermat cost is constructed using $Q_n\subset \Mani$, a sample of $n\geq 1$ i.i.d. points drawn from a common probability measure $\mu$ on $\Mani$ with density $f$. 
Fix $x,y\in\Mani$.  Recall Definition \ref{def:FermatCostGeneral} of the microscopic Fermat costs related to distance $\dist$ and let $\geod_n=\geod_n(x,y)$ by a $C^\dist_n$-geodesic between $x$ and $y$.
We denote by  
$$K_n=K^{\dist}_n(x,y)\coloneqq \card \geod_n (x,y) - 1$$ 
 its number of jumps and by 
$$
J_n=J^{\dist}_n(x,y) \coloneqq \max_{(\geodi_i,\geodi_{i+1})\subset \geod} \dist(\geodi_i,\geodi_{i+1}).
$$
the largest $\dist$-distance between consecutive points in $\geod_n(x,y)$. 
For all these quantities we will  omit, when clarity permits and demands, the dependency on $(x,y)$ and $\dist$. 
Also we will write e.g. $\upb J_n$, $\upb K_n$, etc; to denote corresponding high-probability upper-bounds i.e., upper-bounds satisfying inequalities of the type $\P(J_n\geq \upb J_n)\leq e^{-\const_0n^{\const_1}}$ for some $\const_0,\const_1>0$ and sufficiently large $n$.\par

	Furthermore in the following the constants $\const_0,\const_1>0$ may vary from line to line, as well as the $n$ that is taken sufficiently large in order
	for all bounds to hold. These however, will not depend on any particular pair $x,y\in\Mani$, but will typically depend on a parameter $\eps>0$, and on the remaining  parameters of the model $\Mani$, $\dist$, $f$, $\delta_n$ and $b_n$.\par

 In this work we only deal with piece-wise smooth curves on a manifold $(\Mani,\Rm)$. The length $\length \equiv \length_{\Rm}$ of such a curve
    $\cpath\colon \R\supset [a,b]\to\Mani$
is given by
\begin{align*}
\length \cpath&\coloneqq \int_{[a,b]} \sqrt{\Rm(\cpath'(t),\cpath'(t))}dt\equiv  \int_{[a,b]} \norm{\cpath'(t)}_{\Rm}dt.
%&\equiv \sup_{\norm{P}\to 0} \sum_{i=1}^{\ell-1} \dist_\Mani(\cpath_i,\cpath_{i+1})\\
%&\equiv \lim_{\norm{P}\to 0} \sum_{i=1}^{\ell-1} \dist_\Mani(\cpath_i,\cpath_{i+1}).
\end{align*}
We also let $\length_{\cpath}\colon [a,b]\to\R$ be the corresponding arc-length function along $\cpath$, i.e.
$$
\length_{\cpath}(t) = \int_a^t \norm{\cpath'(s)}_{\Rm}ds.
$$
We will typically work with unit-speed curves (in the $\Rm$-metric) defined on $[0,\length\cpath]$ for which we have $\norm{\cpath'(t)}_{\Rm}=1$ and  $\length_{\cpath}(t)\equiv t.$\par
We also write 
\begin{align}\label{eq:intAlongCurve}
\int_\cpath h(\length_{\cpath}(z),z) \hspace{2pt}d\cpath
&= \int_{[a,b]} h(\length_\cpath(t),\cpath(t)) \sqrt{\Rm_{\cpath(t)}(\cpath'(t),\cpath'(t))}dt
\end{align}
for the integral of a continuous function  $h\colon [0,\infty)\times\Mani\to\R$ along the curve $\cpath$. (On the left we have made a slight abuse of notation by writing $\length_\cpath(z)$ instead of $\length_{\cpath}(t)$ where $z=\cpath(t)$ for some $t\in[a,b]$.) When $\cpath$ is unit-speed, the above integral can be written as 
$$
\int_\cpath h(\length_{\cpath}(z),z) \hspace{2pt}d\cpath = \int_{0}^{\length\cpath} h(t,\cpath(t))dt.
$$

\subsection{Preliminaries concerning $\Mani$ and the sampling on $\Mani$}
The first lemma 
%of this subsection %is proved in Appendix \ref{sec:proofs_conddensityf_DeltanBound_JnAsmp_preliminaryUpperBoundSFD}.and 
is essentially a direct consequence of Assumptions  \ref{cond:continuousMetricSpace} and \ref{cond:boundeddensity}, see \cite{Chavel2006} Chapter IX. 

\begin{lemma}\label{le:conddensityf}
For the pair $(\Mani,\mu)$, there exists a constant $\const_d$ such that
\begin{equation}\label{eq:density}
\quad \forall x\in S,\quad \mu(\ball_{\dist_\Mani}(x,\delta)) = \const_df(x) \delta^d(1 + \lo(1))
\end{equation}
as $\delta\to0$, where the error term $\lo(1)$ is uniform on $x\in\Mani$.
\end{lemma}
Following \cite{aaron2017} (Theorem 1), and using that $\dist_S$ and $\dist_E$ are locally equivalent (see Proposition \ref{prop:euclideaVSgeodesic} later on), we obtain the following for the spacing $\Delta_n$ defined in \eqref{eq:defSpacing}.
\begin{lemma}\label{le:DeltanBound}
    For any $\eps$ there exist  constants $\const_0,\const_1>0$ depending on $f$ and $\eps$ such that, for sufficiently large $n\geq 1$,
    $$
     \P\left( \Delta_n \geq n^{-\frac{1-\eps}{d}}\right) \leq e^{-\const_0 n^{\const_1}}.
    $$
\end{lemma}
\begin{remark}
Recalling Assumption \ref{cond:DeltanBound} on $\delta_n$, Lemma \ref{le:DeltanBound} implies that, for any $\eps>0$ and any constant $\const>0$, there exist  constants $\const_0,\const_1>0$ such that for sufficiently large $n\geq 1$,
$$
\P\left(\frac{\Delta_n}{\const\delta_n}>\eps \right)\leq  e^{-\const_0 n^{\const_1}}.
$$
\end{remark}

\subsection{Preliminaries concerning the microscopic Fermat paths}

The following lemma is proved in Appendix \ref{sec:proofs_conddensityf_DeltanBound_JnAsmp_preliminaryUpperBoundSFD}.
\begin{lemma}\label{le:JnAsmp}
    There exists a number  $$\theta\in(0,1/2)$$ depending only on $\alpha$ and $\maxv A$, such that the maximum step size within a geodesic satisfies
	\begin{equation}\label{eq:JnBound}
		\forall x,y\in\Mani,\quad J^\dist_n(x,y) \asleq \Delta_{n}/\theta. 
	\end{equation}
\end{lemma}
\begin{remark}
	In particular we have
	$
	C^{\dist_\Mani}_n \aseq  C^{\dist_\Mani*}_n
	$
	where $C^{\dist_\Mani*}_n$ is defined as $C^{\dist_\Mani}_n$ in \eqref{eq:defDFermatMicroGeneral} except that the infimum is taken over all $Q_n$-paths $\path$ 
	satisfying $x\qsim y$ and $\max_{1\leq i \leq \card\path} \dist_\Mani(\pathi_i,\pathi_{i-1})\leq \Delta_{n}/\theta$.
\end{remark}

We also provide an upper bound for $K_n$ which is a key technical element for the arguments in the following sections. First, however, we need the following lemma which provides an upper-bound for $C^{\dist_\Mani}_n$ and which is proved in Appendix \ref{sec:proofs_conddensityf_DeltanBound_JnAsmp_preliminaryUpperBoundSFD}. 
\begin{lemma}\label{le:preliminaryUpperBoundSFD}
Let $0<\eps<1$. On the set $\{\Delta_n/\theta \leq 1\}$ we have
\begin{equation*}
   \forall x,y\in\Mani,\quad  C^{\dist_\Mani}_n(x,y)\leq D^{\dist_\Mani}_{n,\alpha}(x,y).%\leq  23^{\alpha+\frac{\minv A}{\log n}}n^{\frac{1-\alpha(1-\eps)}{d}}\dist_\Mani(x,y).
     %\sum_{i=1}^{\lceil n^{1/d}\dist_\Mani(x,y)\rceil}(n^{-\frac{1}{d}}+2\upb\Delta_n)^{\alpha_i}.
\end{equation*}
In particular, under Assumption \ref{cond:DeltanBound}, we have, for any $\eps>0$ and for  constant $\const_0$ depending on $\eps$,
and  constants $\const_1$ depending only on $d,f,\alpha$, such that for sufficiently large $n$,
\begin{align}\label{eq:uppberBoundCn}
&\P\left(\sup_{\substack{\forall x,y\in\Mani\colon \\
                  \dist_\Mani(x,y)\geq \delta_n}}  {\frac{C^{\dist_\Mani}_n(x,y)}{n^{-\frac{\alpha-1}{d}} D_{\infty,\alpha}(x,y)} } \geq \mu(d,\alpha) +\eps\right)\leq e^{-\const_0 (n\delta_n^d)^{\const_1}}.
\end{align}
\end{lemma}

The next result provides a high-probability upper-bound for $K^{\dist_\Mani}_n(x,y)$ that,
importantly,
scales with $D_{\infty,\alpha}(x,y)$. It is of interest in itself in the context of Fermat distances, providing a stronger version of
Lemma 4 in \cite{HwangDamelinHero2016}. 
\begin{proposition}\label{prop:KnAsmp}
For every $x,y\in\Mani$ we have
\begin{equation}%\label{eq:KlowerBound}
K^{\dist_\Mani}_n(x,y) \asgeq \frac{\dist_\Mani(x,y)\theta}{\Delta_n}.
\end{equation}
%where the constant $\theta$ is as in Lemma \ref{le:JnAsmp}.\par
Also, let $\const_K>0$ be the constant defined in Lemma \ref{le:HDHLemma11} below. Define
$$
\upb K_n(x,y)\coloneqq \lfloor\const_{K} 2\mu(d,\alpha)D_{\infty,\alpha}(x,y)n^{\frac{1}{d}}\rfloor.
$$  
Then, under Assumption \ref{cond:DeltanBound}, and for some constants $\const_0>0,\const_1>0$,
and  sufficiently large $n\geq n_0$,
$$
    \P\left(\sup_{\substack{x,y\in\Mani\\\dist_\Mani(x,y)\geq \delta_n}} \frac{K^{\dist_\Mani}_n(x,y)}{\upb K_n(x,y)} \geq 1
    \right)  \leq  e^{-\const_0 n^{\const_1}}.
    $$
%\ahw{   
Furthermore, the constants $\const_0>0,\const_1>0$ and $n_0$ can be chosen
so that the above bound holds uniformly on the class of input functions $A$ satisfying $\maxv A\leq a$ for some $a>0$.
%}
\end{proposition}

The proof of this proposition is provided in Appendix \ref{sec:proof_prop:KnAsmp}.

\subsection{Local properties of macroscopic Fermat optimal paths}
\label{subsec:locpropcinfty}
%\blue{
Let $\kappa>0$
and consider the set of curves 
\begin{align*}
\mathcal F_\kappa \coloneqq \{\cpath\colon  \cpath &\text{ is a continuous and piece-wise smooth curve }
\\ &\text{satisfying } \norm{\nabla_{\cpath'}\cpath'(t)}_{\Rm}\leq \kappa \text{ for a.e. } t \text{ in its domain}\},
\end{align*}
where $\nabla$ is the Levi-Civita connection of $\Rm$.
Since $\cpath$ is continuous and piece-wise smooth, it is continuous and piece-wise absolutely continuous, so absolutely continuous as well.  
In this section we show that the macroscopic $D_{\infty,\alpha}$-optimal curves are contained in $\mathcal F_\kappa$
for a fixed value of
$\kappa=\kappa^*$ that may depend on the parameters of the model $\Mani,f,\alpha$ but that is uniform on $x,y\in\Mani$.
The latter will be used in Proposition \ref{prop:restrictionToBoundedCurvature} below
to show that the infimum in the definition of $C_\infty$ \eqref{eq:macroFermatCost} can be restricted to the set curves
$\mathcal F_{\kappa^*}$. This, in turn, will allow us to compute a metric derivative of the arc-length function uniformly 
along all $C_\infty$-optimal curves, as in Lemma \ref{le:metricDerivativeOfCurves}. The latter will be used in Section \ref{sec:macroscopicLimit} where the mesoscopic approximation of $C_\infty$ is proved.\par
The structure of the proof of the following proposition follows standard arguments in differential geometry and the calculus of variations in Riemannian manifolds; we include them for completion.
\begin{proposition}[Bounded curvature of $D_{\infty,\alpha}$-geodesics]
\label{prop:boundedCurvatureDinftyGeod}
Set
$$
\kappa^* := \frac{\alpha-1}{d}\,\sup_{x\in S}\norm{\specialm{grad}_{\Rm} \ln f(x)}_{\Rm} ,
$$
which is finite by Assumption \ref{cond:boundeddensity}. Then for every pair $x,y\in S$ there exists a
curve attaining the infimum in the definition of $D_{\infty,\alpha}(x,y)$, and any
such optimal curve $\Gamma^*$, parametrized by unit speed, is smooth
and satisfies
$$
\text{for a.e. } t\in[0,\length{\Gamma^*}],\quad \norm{\nabla_{\Gamma^{*\prime}}\Gamma^{*\prime}(t)}_{\Rm} \le \kappa^*.
$$
In particular $\Gamma^*\in\mathcal{F}_{\kappa^*}$. 
%Note that $\kappa^*$ does not
%depend on $x,y$ nor on the function $A$.
\end{proposition}

\begin{proof}
Write $\varphi \coloneqq \ln f^{-\frac{\alpha-1}{d}}$, so that for any piece-wise smooth curve
$\cpath$ (not necessarily of unit speed) we consider the functional
$$
F(\cpath)\coloneqq%\int_0^{\length\cpath}\varphi(\cpath(t))\,dt
\int e^{\varphi(\cpath(s))}\norm{\cpath'(s)}_{\Rm}\,ds .
$$
The right-hand side is the length of $\cpath$ in the conformal metric
$
\tilde \Rm \coloneqq e^{2\varphi}\, \Rm = f^{-\frac{2(\alpha-1)}{d}}\, \Rm .
$
Hence $D_{\infty,\alpha}$ is precisely the geodesic distance of $(S,\tilde \Rm)$.
By Assumption \ref{cond:boundeddensity} the conformal factor $e^{2\varphi}$ is smooth and bounded away from
$0$ and $\infty$, so $\tilde \Rm$ is a smooth Riemannian metric on the compact
manifold $S$. By the Hopf--Rinow theorem (e.g. Theorem 6.19 \cite{Lee2018}) any two points $x,y\in S$ are joined by a
minimizing $\tilde \Rm$-geodesic, and every such minimizer is smooth. This proves the
existence claim, and it only remains to bound the $\Rm$-curvature of these curves.\par
Let $\cpath^*$ be a minimizing $\tilde \Rm$-geodesic, re-parametrized by unit
$\Rm$-speed, $\norm{\cpath^{*\prime}}_{\Rm}\equiv 1$. 
Such a geodesic is a critical point
of the functional $F$, i.e. the first variation $\frac{d}{ds}F(\cpath^*)\vert_{s=0}$ must vanish for any proper
 variation, say $\mathbf \Gamma\colon (-a,a)\times [0,\length\cpath^*]\to\Mani$, of $\cpath^*$ (Theorem 6.3 \cite{Lee2018}). The path $\cpath^*$ being smooth, the integrand in the
definition of $F$ is also smooth; then, the domain of integration being compact we may differentiate under
the integral sign obtaining
$$
\int_0^{\length\cpath^*}dt\, \frac{\partial}{\partial s}\bigg\vert_{s=0}\left(e^{\varphi(\mathbf \Gamma)}\norm{\frac{\partial}{\partial t}\mathbf \Gamma}_{\Rm}\right).
$$
The derivative above can be computed using the product rule. We have, for the first term,
$\frac{\partial}{\partial s}e^{\varphi(\mathbf \Gamma)}=e^{\varphi(\mathbf \Gamma)}\langle \specialm{grad} \varphi, \pder s \mathbf \Gamma \rangle$ where $\nabla$ and $\specialm{grad}$
are the Levi-Civita connection and the gradient associated to $\Rm$. For the second term
a standard computation (see e.g. the proof of Theorem 6.3 \cite{Lee2018}) yields, for a unit-speed smooth path $\cpath^*$, 
$\pder s \vert_{s=0}\norm{\pder t \mathbf \Gamma}_{\Rm}=\langle D_t V,\Gamma^{*\prime}\rangle=
\frac{d}{dt}\langle V, \cpath^{*\prime}\rangle - \langle V, D_t\cpath^{*\prime}\rangle$, where $V(t)=\pder s \mathbf \Gamma(0,t)$
is the variation field of $\mathbf \Gamma$. Combining both terms we have
\begin{align*}
0&=\int_0^{\length\cpath^*} dt\,\frac{\partial}{\partial s}\bigg\vert_{s=0}\left(e^{\varphi(\mathbf \Gamma)}\norm{\frac{\partial}{\partial t}\mathbf \Gamma}_{\Rm}\right) 
\\&= \int_0^{\length\cpath^*} dt\, \left\{e^{\varphi(\cpath^{*})}\langle \specialm{grad} \varphi(\cpath^{*\prime}), V \rangle
 + e^{\varphi(\cpath^*)}\frac{d}{dt}\langle V, \cpath^{*\prime}\rangle - e^{\varphi(\cpath^*)}\langle V, D_t\cpath^{*\prime}\rangle\right\}.
\end{align*}
Applying integration by parts to the middle term, and noting that the boundary term $e^{\varphi(\cpath^*)} \langle V, \Gamma^{*\prime} \rangle\vert_a^b$ vanishes for proper variations, we obtain
$
    \int_a^b e^{\varphi(\cpath^*)}\frac{d}{dt} \langle V, \Gamma^{*\prime} \rangle dt = -\int_a^be^{\varphi(\cpath^*)}\langle \specialm{grad} \varphi(\cpath^{*\prime}), \Gamma^{*\prime} \rangle \langle V, \Gamma^{*\prime} \rangle dt.
$
Substituting above and using the fundamental lemma of the calculus of variations we obtain
$$
D_t\cpath^{*\prime}= \specialm{grad} \varphi(\cpath^{*\prime}) -  \langle \specialm{grad} \varphi(\cpath^{*\prime}), \Gamma^{*\prime} \rangle\Gamma^{*\prime}.
$$
Taking
$\Rm$-norms gives 
$$
\norm{\nabla_{\Gamma^{*\prime}}\Gamma^{*\prime}(t)}_{\Rm}
\leq \norm{\specialm{grad}\varphi (\Gamma^*(t))}_{\Rm}
\leq \kappa^* ,
$$
for every $t$, which is the stated bound.
\end{proof}

\begin{lemma}\label{le:metricDerivativeOfCurves}
The following holds under Assumption \ref{cond:continuousMetricSpace}. Let $\kappa>0$.
For any unit-speed curve $\mathcal F_\kappa\ni \cpath\colon I\to\Mani$, and any $s, t \in I$, the Riemannian distance satisfies:
\begin{equation}\label{eq:geodLengthSmallApprox}
\frac{\dist_{\Mani}(\cpath(t), \cpath(s))}{\length\cpath\vert_{[s,t]}}= \frac{\dist_{\Mani}(\cpath(t), \cpath(s))}{\abs{t-s}} = 1 + \bO(\abs{t-s})
\end{equation}
where the error term $\bO(|t-s|)$ is uniform over the entire class $\mathcal{F}_\kappa$.
\end{lemma}
\begin{proof}
    Let $\cpath$ be an arbitrary curve in $\mathcal F_\kappa$. Since $\dist_\Mani(\cpath(s),\cpath(t))$ is the length of a $\dist_\Mani$-geodesic curve, we immediately have $\length \restr{\cpath}{[s,t]}=\abs{t-s}> \dist_\Mani(\cpath(s),\cpath(t))$.\par
    To obtain a lower bound for the ratio $\frac{\dist_{\Mani}(\cpath(t), \cpath(s))}{|t-s|}$ we first derive an upper-bound for the ambient acceleration $\cpath''$.
    According to the Gauss Formula (e.g. Theorem 8.2 \cite{Lee2018}), the acceleration of a curve on an embedded manifold $\Mani$ can be decomposed into a tangential and an orthogonal component to $\Mani$ as  $$\cpath''(u) = \nabla_{\cpath'}\cpath'(u) + \specialm{II}(\cpath'(u), \cpath'(u))$$
    where $\nabla_{\cpath'}\cpath'$ is the intrinsic geodesic curvature, and $\specialm{II}$ is the Second Fundamental Form of the embedding.
    Taking the Euclidean norm and recalling Assumption \ref{cond:continuousMetricSpace} we obtain 
    $$\norm{\cpath''(u)}_{\R^D} \leq \norm{\nabla_{\cpath'}\cpath'(u)}_{\Rm} + \norm{II(\cpath'(u), \cpath'(u))}_{\R^D}\leq \kappa + \const_{\specialm {II}}.$$
    The latter gives
    $$
    \norm{\cpath'(t)-\cpath'(s)}_{\R^D} \leq \int_s^t \norm{\cpath''(u)}_{\R^D} du \leq ( \kappa + \const_{\specialm {II}})\abs{t-s}.
    $$
    Using that $\cpath'(t)$ and $\cpath'(s)$ are tangent to $M$ and satisfy $\norm{\cpath'(t)}_{\Rm}=1=\norm{\cpath'(s)}_{ \Rm}$ we obtain, from the polarization identity:
    \begin{align*}\langle \cpath'(t),\cpath'(s) \rangle_{\R^D} 
    &= 
    \frac{\norm{\cpath'(t)}_{\Rm}^{2}+\norm{\cpath'(s)}_{\Rm}^{2}}{2}
     - \frac{1}{2}\norm{\cpath'(t)-\cpath'(s)}^2_{\R^D} 
     \\&\geq 1 - \frac{1}{2} ( \kappa + \const_{\specialm {II}})^2\abs{t-s}^2
     \\&= 1- \tilde\kappa \abs{t-s}^2.
     \end{align*}
    Also, the Cauchy-Schwarz inequality gives
    $\norm{\Gamma(t) - \Gamma(s)}_{\R^D} = \norm{\int_{s}^{t} \Gamma'(u) \, du }_{\R^D} \geq \left\langle \int_{s}^{t} \Gamma'(u) \, du, \Gamma'(s) \right\rangle_{\R^D}$; so that, using the bound above, we obtain 
    $$
\norm{\Gamma(t) - \Gamma(s)}_{\R^D} \geq \int_{s}^{t} \left( 1 - \tilde{k} (u-s)^2 \right) du = |t-s| - \frac{ \tilde{k}}{3} |t-s|^3.
    $$
    Because the Riemannian distance is bounded below by the ambient Euclidean distance we conclude 
    $\dist_\Mani(\Gamma(t), \Gamma(s)) \geq |t-s| - \frac{\tilde{k}}{3} |t-s|^3$. This concludes the proof.\par
    \end{proof}
 The following Proposition states that the cost $C_\infty(x,y)$ can be computed even when restricting the infimum appearing in  \eqref{eq:macroFermatCost}
 to the class of curves $\cpath\in\mathcal F_{\kappa^*}$
 whose length is also bounded by $\length\cpath\leq\const_\length\dist_\Mani(x,y)$ for some universal constant $\const_\length>0$. 

\begin{proposition}\label{prop:restrictionToBoundedCurvature}
    We have
    \begin{equation}\label{eq:CinftyWithBoundedCurvature_lowerBound}
    \forall x,y\in\Mani,\quad
        C_\infty(x,y)=\inf_{\substack{\cpath\in\mathcal F_{\kappa^*}\colon\\ x\overset{\cpath}{\sim}y\\ \length\cpath\leq \const_\length\dist_\Mani(x,y)}}\int_0^{\length\cpath} e^{-\frac{A(t)}{d}} f^{-\frac{\alpha-1}{d}}(\cpath(t)) dt
    \end{equation}
    where $\const_\length>0$ is a sufficiently large but fixed constant that is universal for all $x,y\in\Mani$, and the infimum is taken over all unit-speed piece-wise smooth curves $\cpath\in\mathcal F_{\kappa*}$ that connect $x$ and $y$.
    %; and
    %where the error term $\lo(1)$ holds as $n\to\infty$.
\end{proposition}
\begin{proof}
    Fix $x,y\in\Mani$. 
     Let $\eps>0$ be arbitrary and consider $\cgeod\equiv \cpath^*_{x,y}$ a piece-wise smooth curve $\cgeod\colon[0,\length\cgeod]\to\Mani$ satisfying
    $$
    \abs{\int_0^{\length\cgeod} e^{-\frac{A(t)}{d}} f(\cgeod(t))^{-\frac{\alpha-1}{d}}dt-C_\infty(x,y)} \leq \eps,
    $$
    where for simplicity we have assumed that $\cgeod$ is parametrized with unit-speed.\par
    We first prove that the infimum in \eqref{eq:macroFermatCost} may restricted to curves with a bounded length
    \begin{equation}\label{eq:lengthBoundCinfty}
    \length\cgeod(x,y)\leq \const_\length\dist_\Mani(x,y)
    \end{equation}
    for some constant $\const_\length$.
    Indeed, note that
    $$
    e^{-\frac{\maxv A}{d}}(\maxv f)^{-\frac{\alpha-1}{d}}\length\cgeod \leq  C_\infty(x,y) +\eps \leq (\minv f)^{-\frac{\alpha-1}{d}}\dist_\Mani(x,y) + \eps
    $$
where the second inequality holds since $C_\infty(x,y)$ is in particular smaller than the cost of traversing the $\dist_\Mani$-geodesic between $x$ and $y$. Taking $\eps<\dist_\Mani(x,y)$ yields \eqref{eq:lengthBoundCinfty}.\par
    We now prove that we may further restrict the infimum in  \eqref{eq:macroFermatCost} to the set of curves $\mathcal F_{\kappa^*}$ that also satisfy \eqref{eq:lengthBoundCinfty}. Given that $C_\infty(x,y)$ is defined as an infimum over a larger set of curves, we need only prove that that for any $\eps>0$ there exists  a path $\cpath\in\mathcal F_{\kappa^*}$ for which
    \begin{equation}\label{eq:boundedCurvature1}
    \int_0^{\length\cpath} e^{-\frac{A(t)}{d}} f^{-\frac{\alpha-1}{d}}(\cpath(t)) dt \leq (1+\eps) \int_0^{\length\cgeod} e^{-\frac{A(t)}{d}} f(\cgeod(t))^{-\frac{\alpha-1}{d}}dt .
    \end{equation}
    Let $\delta>0$. Define a  partition 
    $\partition=(\cgeod(t_i))_{i=1}^{\ell}$ of $\cgeod$ recursively as follows: let $t_1=0$ so that $\cgeod(t_1)=x$. 
    Assume that $t_i$ has been defined and set
    $$
    t_{i+1} = \inf\{s\geq t_i\colon \dist_\Mani(\cgeod(s),\cgeod(t_i))\geq \delta\}.
    $$
    The sequence of $t_i$'s ends at $t_\ell$ where $\ell\equiv \ell(x,y)$ is an integer-valued function of $(x,y)$. The 
    index $\ell$ is the first time, in the above recursion, that $\dist_\Mani(\cgeod(t_\ell),y)\leq \delta.$
    Note that $\cgeod$ being continuous, we have
    $$
    \forall 1\leq i\leq \ell-1\colon\quad \dist_\Mani(\cgeod(t_i),\cgeod(t_{i+1}))=\delta.
    $$
    Let  $\cpath^*_i\coloneqq \restr{\cpath_{x,y}^*}{[t_i,t_{i+1}]}$. We compute
    \begin{align}\label{eq:CinftyLowerBound1_1}
    \int_0^{\length\cgeod} e^{-\frac{A(t)}{d}} f(\cgeod(t))^{-\frac{\alpha-1}{d}}dt
    &= \sum_{i=1}^{\ell}  \int_{\cpath^*_i} e^{-\frac{A\left(\sum_{j=1}^{i-1} \length\cpath^*_j + \length_{\cpath^*_i}(z)\right)}{d}}f(z)^{-\frac{\alpha-1}{d}} d\cpath^*_i
    \nonumber\\&
    \geq \sum_{i=1}^{\ell-1} e^{-\frac{A\left(\sum_{j=1}^{i} \length\cpath^*_j \right)}{d}}   
     \int_{\cpath^*_i} f(z)^{-\frac{\alpha-1}{d}} d\cpath^*_i
     \nonumber\\&
     \geq \sum_{i=1}^{\ell-1} e^{-\frac{A\left(\sum_{j=1}^{i} \length\cpath^*_j \right)}{d}}   
     D_{\infty,\alpha}(\cgeod(t_{i}),\cgeod(t_{i+1})).
    \end{align}
    Since $\cgeod$ has a curvature bounded by some $\kappa>0$,
    %\blue{letting $\kappa=\kappa_{x,y}$ be a bound on the curvature of $\cgeod(x,y)$ (which always exists since $\cgeod(x,y)$ is assumed to be piece-wise smooth)} we obtain, 
    from eq. \eqref{eq:geodLengthSmallApprox} we have 
    $$
    \forall 1\leq i\leq \ell-1\colon\quad \sum_{j=1}^{i} \length\cpath^*_j = \sum_{j=1}^{i} \dist_\Mani(\cgeod(t_{i}),\cgeod(t_{i+1}))(1+\lo(1))
    $$
    where the error term $\lo(1)$ holds as $\delta\to 0$ and is uniformly on the class of curves $\mathcal F_{\kappa}$.\par
    Similarly, letting $\cpath_i$ be the geodesic curve attaining $D_{\infty,\alpha}(\cgeod(t_i),\cgeod(t_{i+1}))$
    we have, by Proposition \ref{prop:boundedCurvatureDinftyGeod} and Lemma \ref{le:metricDerivativeOfCurves},
    $$
    \forall 1\leq i\leq \ell-1\colon\quad \sum_{j=1}^{i} \length\cpath_j = \sum_{j=1}^{i} \dist_\Mani(\cgeod(t_{i}),\cgeod(t_{i+1}))(1+\lo(1))
    $$
    where the error term $\lo(1)$ holds as $\delta\to 0$ and is uniform on the class of curves $\mathcal F_{\kappa^*}.$
    Let $\cpath$ be the concatenation $\cpath=\cpath_1+\cdots+\cpath_{\ell-1}$. The latter implies
    \begin{align*}
    &\length\cpath
    = \sum_{i=1}^{\ell-1} \dist_\Mani(\cgeod(t_{i}),\cgeod(t_{i+1}))(1+\lo(1))
    \\&\leq  \sum_{i=1}^{\ell-1} \length\cpath^*_j (1+\lo(1))
    \leq  \length\cgeod (1+\lo(1))
    \end{align*}
    so that $\length\cpath\leq \const_\length\dist(x,y)$ for some sufficiently large universal constant $\const_\length>0$.
    From the latter we furthermore  obtain
    $$
    \max_{1\leq i\leq \ell-1}\left\{\sum_{j=1}^{i} \length\cpath^*_j 
    - \sum_{j=1}^{i} \length\cpath_j\right\}\leq (2\length\cgeod)\lo(1).  
    $$
    In particular
    $$
    \max_{1\leq i\leq \ell-1}\left\{A\left(\sum_{j=1}^{i} \length\cpath^*_j 
    \right)- A\left(\sum_{j=1}^{i} \length\cpath_j\right)\right\}=\modulus_A(\lo(1)) 
    $$
    where the error term $\lo(1)$ holds as $\delta\to0$ and is uniform on $\mathcal F_\kappa \cup \mathcal F_{\kappa^*}$. 
    Plugging the above in \eqref{eq:CinftyLowerBound1_1} we obtain
    \begin{equation}\label{eq:boundedCurvature2}
    \int_0^{\length\cgeod} e^{-\frac{A(t)}{d}} f(\cgeod(t))^{-\frac{\alpha-1}{d}}dt
    \geq e^{-\frac{\modulus_A(\lo(1))}{d}} \sum_{i=1}^{\ell-1} e^{-\frac{A\left(\sum_{j=1}^{i} \length\cpath_j \right)}{d}}   
     \int_{\cpath_i} f(z)^{-\frac{\alpha-1}{d}} d\cpath_i.
    \end{equation}
    An analogous argument also yields
    \begin{equation}\label{eq:boundedCurvature3}
    \sum_{i=1}^{\ell-1} e^{-\frac{A\left(\sum_{j=1}^{i} \length\cpath_j \right)}{d}}   
     \int_{\cpath_i} f(z)^{-\frac{\alpha-1}{d}} d\cpath_i
    \geq e^{-\frac{\modulus_A(\lo(1))}{d}}\int_\cpath e^{-\frac{A(\length_{\cpath}(z))}{d}} f(\cpath(z))^{-\frac{\alpha-1}{d}} d\cpath
    \end{equation}
    where still the error term $\lo(1)$ holds as $\delta\to0$ and is uniform on $\mathcal F_\kappa \cup \mathcal F_{\kappa^*}$. Putting together \eqref{eq:boundedCurvature2} and \eqref{eq:boundedCurvature3}, and recalling  Assumption \ref{cond:functionA},  we obtain \eqref{eq:boundedCurvature1}.
\end{proof}
%}

\section{Step 1: Microscopic to mesoscopic}\label{sec:mesoscopicApproximation}

In this section we use the notation $C_n=C_n^{\dist_\Mani}$ while $K_n$ will refer to $K^{\dist_\Mani}_n$.
For fixed $x,y\in\Mani$, consider the following set of sequences of points  
\begin{align}\label{eq:defUn}
    \mathcal U_n(x,y) \coloneqq 
    \Big\{
    (u_k)_{k=1}^\ell\subset \Mani \colon &\ell\leq n +2,\nonumber\\
    & u_1=x, \quad \dist_\Mani(y,u_\ell)\le{\delta_n} ;
    \nonumber\\
    &\forall 2\leq k\leq\ell,\quad \delta_n \leq \dist_\Mani(u_{k-1},u_{k})%\leq {\delta_n+4\Delta_n/\theta}%\leq \delta_n + \upb\Delta_n/\theta
    \Big\}.
\end{align}
These are mesoscopic discrete paths of points in $\Mani$, starting at $x$, and ending close to $y$. Consecutive points in the path are at a mesoscopic distance of at least $\delta_n$.\par
 We define the following random mesoscopic costs between a pair of points $x,y\in\Mani$:
$$
\widehat C_n(x,y) \coloneqq \inf_{\substack{(u_k)_{k=1}^{\ell} \in \mathcal U_n(x,y)\\ (u_k)_{k=1}^{\ell}\subset Q_n \cup \{x,y\}}}
\sum_{k=1}^{\ell-1} D_{n,\alpha_{n,k}}(u_k,u_{k+1}),
$$
and 
$$
\widehat C_{n,\infty}(x,y) \coloneqq \inf_{\substack{(u_k)_{k=1}^{\ell} \in \mathcal U_n(x,y)\\ (u_k)_{k=1}^{\ell}\subset Q_n \cup \{x,y\}}}
\sum_{k=1}^{\ell-1} e^{-\frac{A(k\delta_n)}{d}}D_{\infty,\alpha}(u_k,u_{k+1}).
$$

The objective of this section is to prove the following result, which is the first main step of the proof of our main theorem, Theorem \ref{th:FermatCostConvergence}.
\begin{theorem}\label{theo:hatCnVshatCninfty}
There exist constants $\const_0>0$ and $n_0\geq0$ that may depend on $\eps$, 
    and
    $\const_1>0$ depending only on  $f,d,\alpha,\maxv A$ such that, under Assumption \ref{cond:bn}, %\ref{cond:DeltanBound} 
$$
\forall n\geq n_0,\quad
\P\left( \sup_{\substack{x,y\in\Mani\\\dist_\Mani(x,y)\geq b_n}}
\abs{\frac{n^{\frac{\alpha-1}{d}} C_n(x,y)}
     {\widehat C_{n,\infty}(x,y) }
     -\mu(d,\alpha)}\ge \eps
\right) \leq e^{-\const_0 n^{\const_1}}.
$$
\end{theorem}

\subsection{Shape Theorem for $D_{n,\alpha+\frac{\bar \alpha}{\log n}}$}
The following result is an extension of Theorem 1 \cite{HwangDamelinHero2016} to the case where the parameter $\alpha$ therein scales with $n$.
This is a key result used in section \ref{sec:mesoscopicApproximation} below to translate the microscopic percolation mechanism in \ref{eq:defDFermatMicroGeneral} into
a mesoscopic one, where now the cost of a given step within a path is of the form $e^{-\frac{A(k\delta_n)}{d}}D_{\infty,\alpha}(u_k,u_{k+1})$.

 Theorem 1 of \cite{HwangDamelinHero2016} (see also Theorem 6 of \cite{FernandezBorghiniMindlinGroisman2023}) states that, under Assumptions \ref{cond:continuousMetricSpace}, \ref{cond:DeltanBound} and \ref{cond:boundeddensity} we have for sufficiently large $n$,
    \begin{equation}\label{eq:damelin_concentration}
 \P\left(\sup_{\substack{\forall x,y\in\Mani\colon \\
                  \dist_\Mani(x,y)\geq \delta_n}}  \left|{\frac{D_{n,\alpha}(x,y)}{n^{-\frac{\alpha-1}{d}} D_{\infty,\alpha}(x,y)} } -\mu(d,\alpha)\right|\geq \eps \right)\leq e^{-\const_0 (n\delta_n^d)^\alpha},
    \end{equation}
    where $\const_0$ depends on $\eps$. The next proposition gives a similar result when instead of a fixed parameter $\alpha$, we consider $\alpha+\frac{\bar \alpha}{\log n}$ which varies with $n$.

\begin{proposition}\label{prop:alphanDnConv}
    Let $\alpha>1$ and $0\leq\bar\alpha\leq\maxv A$. Fix $\eps>0$. Then there exist constants $\const_0>0$ that may depend  on $\eps$, 
    and
    $\const_1>0$  depending only on  $f,d,\alpha,\maxv A$ (but not on $\bar\alpha$) such that, under Assumption \ref{cond:DeltanBound}, we have for sufficiently large $n$,
    $$
     \P\left(\sup_{\substack{x,y\in\Mani\\\dist_\Mani(x,y)\geq \delta_n}}
    \abs{\frac{n^{\frac{\alpha-1}{d}}D_{n,\alpha+\frac{\bar \alpha}{\log n}}(x,y)}{e^{-\frac{\bar\alpha}{d}}D_{\infty,\alpha}(x,y)} - \mu(d,\alpha)} \ge \eps
    \right)
    \leq e^{-\const_0n^{\const_1}}.
    $$
\end{proposition}

\begin{proof}
Fix $n\ge1$, $\alpha>1$, $\bar\alpha\in[0,\maxv A]$, $\eps>0$ and $x,y\in\Mani$. Denote by $\tilde\cpath$ the set of $Q_n$-paths $\path$ 
    satisfying $x\qsim y$ and $\max_{1\leq i \leq \card\path} \dist_\Mani(\pathi_i,\pathi_{i-1})\leq \Delta_{n}/\theta$ where we recall that $\theta<1$. We have
\begin{align*}
D_{n,\alpha+\frac{\bar \alpha}{\log n}}(x,y)&\leq \inf_{\path\in \tilde\Gamma}\sum_{i=1}^{\card \path -1}\dist_\Mani(\pathi_{i+1},\pathi_i)^{\alpha} \dist_\Mani(\pathi_{i+1},\pathi_i)^{\frac{\bar \alpha}{\log(n)}}\\&\leq \theta^{-\frac{\bar \alpha}{\log(n)}}\Delta_n^{\frac{\bar \alpha}{\log(n)}}\inf_{\path\in \tilde\Gamma}\sum_{i=1}^{\card \path -1}\dist_\Mani(\pathi_{i+1},\pathi_i)^{\alpha}.
% \\&= \Delta_n^{\frac{\bar \alpha}{\log(n)}}D_{n,\alpha}(x,y).
\end{align*}
According to Lemma \ref{le:JnAsmp}, the set
\begin{equation}
\mathcal C_n:=\left\{D_{n,\alpha}(x,y)=\inf_{\path\in \tilde\Gamma}\sum_{i=1}^{\card \path -1}\dist_\Mani(\pathi_{i+1},\pathi_i)^{\alpha}\right\}
\end{equation}
has probability 1. 
Let us denote
$$\mathcal E_n:=\left\{\sup_{\substack{x,y\in\Mani\\\dist_\Mani(x,y)\geq \delta_n}}\frac{n^{\frac{\alpha-1}{d}}D_{n,\alpha+\frac{\bar \alpha}{\log n}}(x,y)}{e^{-\frac{\bar\alpha}{d}}D_{\infty,\alpha}(x,y)}-\mu(d,\alpha)>\epsilon\right\}.$$
For any $\beta>0$, let $\mathcal A_n:=\{\Delta_n\leq n^{-\frac{1-\beta}{d}}\}$, so $\Delta_n^{\frac{\bar \alpha}{\log(n)}}\leq e^{\frac{\beta\maxv A}{d}}e^{-\frac{\bar\alpha}{d}}$ and 
\begin{align*}
&\P(\mathcal E_n\cap\mathcal A_n)=\P(\mathcal E_n\cap\mathcal A_n\cap \mathcal C_n)\\&\leq \P\left(\sup_{\substack{x,y\in\Mani\\\dist_\Mani(x,y)\geq \delta_n}}\frac{n^{\frac{\alpha-1}{d}}D_{n,\alpha}(x,y)}{D_{\infty,\alpha}(x,y)}>\theta^{\frac{\maxv A}{\log(n)}}e^{-\frac{\beta\maxv A}{d}}(\mu(d,\alpha)+\epsilon)\right).
\end{align*}
If $\beta<\frac{d}{\maxv A}\log\left(\frac{\mu(d,\alpha)+\epsilon}{\mu(d,\alpha)}\right)$, we can use \eqref{eq:damelin_concentration} to conclude that for sufficiently large $n\geq n_0$, $\P(\mathcal E_n\cap\mathcal A_n)\leq e^{-\const_0 (n\delta_n^d)^\alpha}$ for some $\const_0>0$ not depending on $\bar \alpha$ (but possibly on $\maxv A$).
%depending only on $\eps,\beta,\alpha,\maxv A$ and $d$. 
According to Lemma \ref{le:DeltanBound}, a similar bound holds  for $\P(\mathcal A_n^c)$. We thus conclude 
\begin{equation}\label{eq:En}
\P(\mathcal E_n)\leq \P(\mathcal A_n^c)+\P(\mathcal E_n\cap\mathcal A_n) \leq  e^{-\const_0 n^{\const_1}}
\end{equation}
for some constants $\const_0,\const_1$ and for sufficiently large $n\geq n_0$ not depending on $\bar \alpha$ (but possibly on $\maxv A$). 

We now work on obtaining an upper bound on the probability of 
$$\mathcal F_n:=\left\{ \inf_{\substack{x,y\in\Mani\\\dist_\Mani(x,y)\geq \delta_n}}\frac{n^{\frac{\alpha-1}{d}}D_{n,\alpha+\frac{\bar \alpha}{\log n}}(x,y)}{e^{-\frac{\bar\alpha}{d}}D_{\infty,\alpha}(x,y)}-\mu(d,\alpha)<-\epsilon\right\}.$$
For any fixed $x,y$, consider a $\path^*=\path^*(x,y)$ reaching the infimum in the definition of $D_{n,\alpha+\frac{\bar \alpha}{\log n}}(x,y)$. We have for any $\beta>0$

\begin{align*}
&D_{n,\alpha+\frac{\bar \alpha}{\log n}}(x,y)
=\sum_{i=1}^{\card \path^* -1}\dist_\Mani(\pathi^*_{i+1},\pathi^*_i)^{\alpha} \dist_\Mani(\pathi^*_{i+1},\pathi^*_i)^{\frac{\bar \alpha}{\log(n)}}\\
%&\geq \sum_{i=1}^{\card \path^* -1}{\bf 1}\{\dist_\Mani(\pathi^*_{i+1},\pathi^*_i)>n^{-\frac{(1+\beta)}{d}}\}\dist_\Mani(\pathi^*_{i+1},\pathi^*_i)^{\alpha} \dist_\Mani(\pathi^*_{i+1},\pathi^*_i)^{\frac{\bar \alpha}{\log(n)}}\\
&\geq n^{-\frac{(1+\beta)}{d}\frac{\bar \alpha}{\log(n)}}\sum_{i=1}^{\card \path^* -1}{\bf 1}\{\dist_\Mani(\pathi^*_{i+1},\pathi^*_i)>n^{-\frac{(1+\beta)}{d}}\}\dist_\Mani(\pathi^*_{i+1},\pathi^*_i)^{\alpha}\\
%&\geq n^{-\frac{(1+\beta)}{d}\frac{\bar \alpha}{\log(n)}}\left(\sum_{i=1}^{\card \path^* -1}\dist_\Mani(\pathi^*_{i+1},\pathi^*_i)^{\alpha}-\sum_{i=1}^{\card \path^* -1}{\bf 1}\{\dist_\Mani(\pathi^*_{i+1},\pathi^*_i)\leq n^{-\frac{(1+\beta)}{d}}\}\dist_\Mani(\pathi^*_{i+1},\pathi^*_i)^{\alpha}\right)\\
&\geq n^{-\frac{(1+\beta)}{d}\frac{\bar \alpha}{\log(n)}}\left(D_{n,\alpha}(x,y)-\sum_{i=1}^{\card \path^* -1}{\bf 1}\{\dist_\Mani(\pathi^*_{i+1},\pathi^*_i)\leq n^{-\frac{(1+\beta)}{d}}\}\dist_\Mani(\pathi^*_{i+1},\pathi^*_i)^{\alpha}\right)\\
&\geq e^{-\frac{\bar\alpha(1+\beta)}{d}}\left(D_{n,\alpha}(x,y)-K_n(x,y) n^{-\frac{\alpha(1+\beta)}{d}}\right).
\end{align*}
On $\mathcal B_n:=\{\sup_{\substack{x,y\in\Mani\\\dist_\Mani(x,y)\geq \delta_n}}\frac{K_n (x,y)}{\upb K_n(x,y)}<1\}$, we have (see Proposition \ref{prop:KnAsmp} for the definition of $\upb K_n$)
$$D_{n,\alpha+\frac{\bar \alpha}{\log n}}(x,y)\geq e^{-\frac{\bar\alpha(1+\beta)}{d}}\left(D_{n,\alpha}(x,y)-\lfloor\const_{K} 2\mu(d,\alpha)D_{\infty,\alpha}(x,y)\rfloor  n^{-\frac{\alpha(1+\beta)-1}{d}}\right)$$
so that $\P(\mathcal F_n\cap \mathcal B_n)$ is bounded by
\begin{align*}
%\P(\mathcal F_n\cap \mathcal B_n)&\leq \P\left(\inf_{\substack{x,y\in\Mani\\\dist_\Mani(x,y)\geq \delta_n}}\frac{n^{\frac{\alpha-1}{d}}D_{n,\alpha}(x,y)}{D_{\infty,\alpha}(x,y)}-{2\const_K \mu(d,\alpha) n^{-\frac{\alpha\beta}{d}}}<e^{\frac{\bar\alpha\beta}{d}}(\mu(d,\alpha)-\epsilon)\right)\\
%&\leq 
\P\left(\inf_{\substack{x,y\in\Mani\\\dist_\Mani(x,y)\geq \delta_n}}\frac{n^{\frac{\alpha-1}{d}}D_{n,\alpha}(x,y)}{D_{\infty,\alpha}(x,y)}<e^{\frac{\bar\alpha\beta}{d}}(\mu(d,\alpha)-\epsilon)+{2\const_K \mu(d,\alpha) n^{-\frac{\alpha\beta}{d}}}\right).
\end{align*}
Observe that ${2\const_K \mu(d,\alpha) n^{-\frac{\alpha\beta}{d}}}$ vanishes as $n$ diverges. So if we take $\beta<\frac{d}{\maxv A}\log\left(\frac{\mu(d,\alpha)}{\mu(d,\alpha)-\eps}\right)$  we  can use \eqref{eq:damelin_concentration} to conclude that for any sufficiently large $n\geq n_0$, $\P(\mathcal F_n\cap \mathcal B_n)\leq e^{-\const_0 (n\delta_n^d)^\alpha}$ where the constant $\const_0>0$ 
is valid for any $\bar\alpha\leq\maxv A$. Similarly, using Proposition \ref{prop:KnAsmp}, there exist constants $\const_0,\const_1>0$  that are valid for any $\bar\alpha\leq\maxv A$ such that  $\P(\mathcal B_n^c)\leq e^{-\const_0 n^{\const_1}}$ for sufficiently large $n$. This gives
$\P(\mathcal F_n)\leq \P(\mathcal B_n^c)+\P(\mathcal F_n\cap \mathcal B_n)\leq e^{-\const_0 n^{\const_1}}.$ The latter, together with \eqref{eq:En}, give the Proposition. 
\end{proof}

\subsection{From Microscopic to Mesoscopic Percolation}
We have the following.
\begin{lemma}\label{le:CnvshatCn}
	Fix $\eps>0$. Then there exist $\const_0>0,\const_1>0$ such that, for sufficiently large $n$,
	$$
	\P\left(\sup_{\substack{x,y\in\Mani\\\dist_\Mani(x,y)\geq b_n}}\left|\frac{C_n(x,y)}{\widehat C_n(x,y)}-1\right| > \eps \right) \leq e^{- \const_0n^{\const_1}}.
	$$
\end{lemma}
\begin{proof}
 By Lemmas \ref{le:DeltanBound} and \ref{le:JnAsmp} we may work on the high-probability event $\{J_n\leq \Delta_n / \theta < 1\}$. On this event, any discrete path that is a $C_n$- or $D_{n,\alpha}$-geodesic has all consecutive points at distance less than $1$.
\par
Let us first prove the lower bound. Fix $x, y \in \Mani$ with $\dist_{\Mani}(x,y) \ge b_n$
and let $(\geod,\iota)$ be the optimal path and marks giving $C_n(x,y)$, with $\ell=\card\iota$. Note that, on the event $\{\Delta_n / \theta < 1\}$, we may assume  
$(\geod_{\iota(1)},\cdots,\geod_{\iota(\ell)})\in\mathcal U_n$ without loss of generality. Indeed, otherwise we need only append
 $\card\path$ to $\iota$, which keeps the marking admissible
 and does not change the associated cost.
Then, since $(\geod_{\iota(1)},\cdots,\geod_{\iota(\ell)})\in\mathcal U_n$ we directly obtain
$$
C_n(x,y)\geq  \sum_{k=1}^{\ell-1} D_{n,\alpha_k}(\geod_{\iota(k)}, \geod_{\iota(k+1)}) \geq  \widehat{C}_n(x,y).
$$
This proves the lower bound.
We now prove the upper bound. 
%For any $\eta > 0$ 
Choose a mesoscopic path $(u_1, \dots, u_\ell) \in \mathcal{U}_n(x,y)$ such that
$$
\sum_{k=1}^{\ell-1} D_{n,\alpha_k}(u_k, u_{k+1}) = \widehat{C}_n(x,y).
$$

For each $k = 1, \dots, \ell-1$, let $\pi_k$ be a $D_{n,\alpha_k}$-geodesic from $u_k$ to $u_{k+1}$. 
Concatenate these paths to obtain a path $\pi$ from $x$ to $u_\ell$, and append the direct step from $u_\ell$ to $y$. Call the resulting path $\path$, and let $\iota(k)$ be defined
as the index $i$  such that $\path_{i}$ is the starting point of the geodesic $\pi_k$ (formally $\iota(k)=1+\sum_{j=1}^{k-1}(\card\pi_j-1)$). Then, from the definition of $\mathcal U_n$, we have $\iota\in M_{\dist,\delta_n}(\path)$. Thus 
$$
C_n(x,y) \leq \sum_{k=1}^{\ell-1} D_{n,\alpha_k}(u_k, u_{k+1}) + \dist(u_\ell, y)^{\alpha_\ell}.
$$
Since $\dist(u_\ell, y) < \delta_n$ and $\alpha_\ell \ge \alpha$, we have $\dist(u_\ell, y)^{\alpha_\ell} \le \delta_n^{\alpha}$ (recall $\delta_n < 1$ for large $n$). Consequently,
$
C_n(x,y) \le \widehat{C}_n(x,y) +  \delta_n^{\alpha}.
$

%As $\eta > 0$ was arbitrary we obtain 
%$
%C_n(x,y) \leq \widehat{C}_n(x,y) + \delta_n^{\alpha}.
%$
\par
Putting the above two both bounds together we have proved that for all $x,y\in\Mani$ with $\dist(x,y) \ge b_n$ we have, on the event $\{\Delta_n / \theta < 1\}$, that
$$
\Bigl| \frac{C_n(x,y)}{\widehat{C}_n(x,y)} - 1 \Bigr| \le \frac{\delta_n^{\alpha}}{\widehat{C}_n(x,y)}.
$$
Let us now prove that  \begin{equation}\label{eq:hatCnVSdeltal_n^alpha}
    \P\left(\sup_{\substack{x,y\in\Mani\\\dist_\Mani(x,y)\geq b_n}} \frac{\delta_n^\alpha}{\widehat C_n(x,y)} > \eps \right) \leq e^{- \const_0 n^{\const_1}}
    \end{equation}
    for large enough $n$. Note that on the event $\{\Delta/\theta < 1\}$ we have $$D_{n,\alpha+\frac{\maxv A}{\log n}}(u_k,u_{k+1})\leq D_{n,\alpha+\frac{A(k\delta_n)}{\log n}}(u_k,u_{k+1})$$ for
	all $u_k,u_{k+1}\in\Mani$. Then, for arbitrary $x,y\in\Mani$ such that $\dist_\Mani(x,y)\geq \delta_n$ we have  $\widehat C_n(x,y)>D_{n,\alpha+\frac{\maxv A}{\log n}}(x,u)$ for some $u$ satisfying $\dist_\Mani(x,u)\geq \dist_\Mani(x,y)-\delta_n.$ 
	Furthermore, using Proposition \ref{prop:alphanDnConv} we obtain that,  for some constant $\const$ and for large enough $n\geq1$,
	\begin{align*}
		&\P\left(\inf_{\substack{x,u\in\Mani\\\dist_\Mani(x,u)\geq \delta_n}} \frac{n^{\frac{1-\alpha}{d}}\const D_{\infty,\alpha}(x,u)}{D_{n,\alpha+\frac{\maxv A}{\log n}}(x,u)}\geq 1\right)
		\\&=\P\left(\sup_{\substack{x,u\in\Mani\\\dist_\Mani(x,u)\geq \delta _n}} \frac{n^{\frac{\alpha-1}{d}}D_{n,\alpha+\frac{\maxv A}{\log n}}(x,u)}{\const D_{\infty,\alpha}(x,u)}\leq 1\right)
		\leq e^{- \const_0 n^{\const_1}}.
	\end{align*}
	Thus, as soon as $\dist_\Mani(x,y)-\delta_n\geq b_n-\delta_n\geq\delta_n$ we obtain, with probability at least $1- e^{- \const_0 n^{\const_1}}$, that
	$$
	\frac{\delta_n^\alpha}{\widehat C_n(x,y)}  \leq \frac{\delta_n^\alpha}{n^{\frac{1-\alpha}{d}}\const D_{\infty,\alpha}(x,u)}
	\leq \frac{\delta_n^\alpha}{n^{\frac{1-\alpha}{d}}\const (\maxv f)^{-\frac{\alpha-1}{d}} b_n}.
	$$
	Taking $b_n$ satisfying $\frac{\delta_n^\alpha}{n^{\frac{1-\alpha}{d}}b_n}\to 0$, in particular under Assumption \ref{cond:bn}, 
	we  obtain \eqref{eq:hatCnVSdeltal_n^alpha}.
\end{proof}

\subsection{Proof of Theorem \ref{theo:hatCnVshatCninfty}}

We are now ready to prove Theorem \ref{theo:hatCnVshatCninfty} as a consequence of Proposition \ref{prop:alphanDnConv} and Lemma \ref{le:CnvshatCn}.
\begin{proof}
	Recall $\alpha_k\equiv\alpha_{n,k}=\alpha+\frac{A(k\delta_n)}{\log n}.$
    We claim that the event
    $$
    \mathcal E_n=\left\{ \sup_{\substack{x,y\in\Mani\\\dist_\Mani(x,y)\geq \delta_n}}
\abs{\frac{n^{\frac{\alpha-1}{d}}\widehat C_n(x,y)}
     {\widehat C_{n,\infty}(x,y) }
     -\mu(d,\alpha)}\ge \eps\right\} % ; \sup_{\substack{x,y\in\Mani\\\dist_\Mani(x,y)\geq \delta_n}}K_n \leq  \upb K_n
    $$
    is contained in the event 
    $$
    \mathcal E_n'=\bigcup_{k=1}^{n+2} \left\{\sup_{\substack{x',y'\in\Mani\colon\\\dist_\Mani(x',y')\geq \delta_n}}
                        \abs{\frac{n^{\frac{\alpha-1}{d}}D_{n,\alpha_k}(x',y')}
                             {e^{-\frac{A(k\delta_n)}{d}}D_{\infty,\alpha}(x',y')} - \mu(d,\alpha)} \geq \eps  
     \right\}.
    $$
(We prove that $(\mathcal E_n')^c\subset \mathcal E_n^c$ at the end of the proof.) Then using Proposition \ref{prop:alphanDnConv} %and Proposition \ref{prop:KnAsmp}
we see that, for some constants $\const_0>0,\const_1>0$ that depend only on $\maxv A$ but not on $k\geq 1$ below, we have
    \begin{align}\label{eq:void}
    &\P\left( \sup_{\substack{x,y\in\Mani\\\dist_\Mani(x,y)\geq \delta_n}}
\abs{\frac{n^{\frac{\alpha-1}{d}}\widehat C_n(x,y)}
     {\widehat C_{n,\infty}(x,y) }
     -\mu(d,\alpha)}\ge \eps
\right)\notag\\
    &\leq \sum_{k=1}^{n+2} \P\left(\sup_{\substack{x',y'\in\Mani\colon\\\dist_\Mani(x',y')\geq \delta_n}}
                        \abs{\frac{n^{\frac{\alpha-1}{d}}D_{n,\alpha_k}(x',y')}
                             {e^{-\frac{A(k\delta_n)}{d}}D_{\infty,\alpha}(x',y')} - \mu(d,\alpha)} \geq \eps \right)
     \nonumber\\
     &\leq       \sum_{k=1}^{n+2} e^{-\const_0n^{\const_1}}            
      = n e^{-\const_0n^{\const_1}}.
    \end{align}
Assuming the above, the proof of Theorem \ref{theo:hatCnVshatCninfty} is completed by writing 
    \begin{align*}
\abs{\frac{n^{\frac{\alpha-1}{d}} C_n(x,y)}
     {\mu(d,\alpha)\widehat C_{n,\infty}(x,y) } -1 } &\leq 
     \abs{\frac{C_n(x,y)}{\widehat C_n(x,y)} - 1}\abs{\frac{n^{\frac{\alpha-1}{d}}\widehat C_n(x,y)}{\mu(d,\alpha)\widehat C_{n,\infty}(x,y)} -1 } \\&+  \abs{\frac{C_n(x,y)}{\widehat C_n(x,y)} - 1} +\abs{\frac{n^{\frac{\alpha-1}{d}}\widehat C_n(x,y)}{\mu(d,\alpha)\widehat C_{n,\infty}(x,y)} -1 }, 
    \end{align*}
    and then using \eqref{eq:void} together with Lemma \ref{le:CnvshatCn} to bound each term with high probability.
It only remains to prove that $(\mathcal E_n')^c\subset \mathcal E_n^c$.    Recall that by the definition
    of $\mathcal U_n$ we have $\ell\leq n+2$ whenever $(u_k)_{k=1}^\ell\in\mathcal U_n$.
    Thus,
    on the event $(\mathcal E_n')^c$ we have, for any collection 
    $(u_k)_{k=1}^\ell\in\mathcal U_n$ and any $1\leq k\leq \ell-1$, that
    \begin{align*}
    &\abs{n^{\frac{\alpha-1}{d}}D_{n,\alpha_k}(u_k,u_{k+1})
                              - \mu(d,\alpha)e^{-\frac{A(k\delta_n)}{d}}D_{\infty,\alpha}(u_k,u_{k+1})} 
                              \\&\leq \eps e^{-\frac{A(k\delta_n)}{d}}D_{\infty,\alpha}(u_k,u_{k+1}).
    \end{align*}
    Therefore $\forall (u_k)_{k=1}^\ell \in\mathcal U_n;$
    \begin{align}\label{eq:hatCn1}
    &\abs{
    \sum_{k=1}^{\ell-1}
    n^{\frac{\alpha-1}{d}}D_{n,\alpha_k}(u_k,u_{k+1})
                              - \mu(d,\alpha)\sum_{k=1}^{\ell-1} e^{-\frac{A(k\delta_n)}{d}}D_{\infty,\alpha}(u_k,u_{k+1})
    }
    \nonumber\\
    &\leq \eps \sum_{k=1}^{\ell-1} e^{-\frac{A(k\delta_n)}{d}}D_{\infty,\alpha}(u_k,u_{k+1}).
    \end{align}\par
    Fix arbitrary $x,y\in\Mani$ satisfying $\dist_\Mani(x,y)\geq\delta_n$. 
    Equation \eqref{eq:hatCn1} implies, for any $(u_k)_{k=1}^\ell\in \mathcal U_n(x,y)$, that  on the event $(\mathcal E_n')^c$ we have
    \begin{align*}
    &(\mu(d,\alpha)-\eps)\sum_{k=1}^{\ell-1} e^{-\frac{A(k\delta_n)}{d}}D_{\infty,\alpha}(u_k,u_{k+1})
    \leq n^{\frac{\alpha-1}{d}}\sum_{k=1}^{\ell-1}
    D_{n,\alpha_k}(u_k,u_{k+1}) \\&\leq (\mu(d,\alpha)+\eps)\sum_{k=1}^{\ell-1} e^{-\frac{A(k\delta_n)}{d}}D_{\infty,\alpha}(u_k,u_{k+1}).
    \end{align*}
    Taking the infimum over $(u_k)_{k=1}^\ell\in\mathcal U_n(x,y)\cap \{Q_n\cup \{x\}\}$ we obtain
    %, on the event
    %$\{\sup_{\substack{x,y\in\Mani\\\dist_\Mani(x,y)\geq \delta_n}}K_n \leq  \upb K_n\}$, that
    \begin{align*}
    (\mu(d,\alpha)-\eps)\widehat C_{n,\infty}(x,y) \leq n^{\frac{\alpha-1}{d}}\widehat C_n(x,y) \leq 
    (\mu(d,\alpha)+\eps)\widehat C_{n,\infty}(x,y);
    \end{align*}
    i.e.
    $
    \abs{\frac{n^{\frac{\alpha-1}{d}}\widehat C_n(x,y)}{\widehat C_{n,\infty}(x,y)}-\mu(d,\alpha)}\leq \eps.
    $
    Since the pair $x,y\in\Mani$ satisfying $\dist_\Mani(x,y)\geq\delta_n$ was arbitrary, we conclude
    $
    \sup_{\substack{x,y\in\Mani\\\dist_\Mani(x,y)\geq \delta_n}}
    \abs{\frac{n^{\frac{\alpha-1}{d}}\widehat C_n(x,y)}{\widehat C_{n,\infty}(x,y)}-\mu(d,\alpha)}\leq \eps.
    $
    This proves $(\mathcal E_n')^c\subset \mathcal E_n^c$.
    The proof of the theorem is therefore concluded.
\end{proof}

\section{Step 2: Mesoscopic to macroscopic}\label{sec:macroscopicLimit}

The  objective of this section is to prove Theorem \ref{th:hatCnInfVSCInf} below, which is the second main step towards our proof of Theorem \ref{th:FermatCostConvergence}. This step can be seen as a quantification of how well the mesoscopic cost $\widehat C_{n,\infty}$ approximates the macroscopic limit $C_\infty$.

\begin{theorem}[$\widehat C_{n,\infty}$ approximates $C_\infty$]\label{th:hatCnInfVSCInf}
The following holds under Assumption \ref{cond:functionA}. For any $\eps>0$ there exist constants $\const_0>0,\const_1>0$ 
    such that for large enough $n$,
    \begin{equation}\label{eq:hatCnInftyVSCinfty}
        \P\left(\sup_{\substack{x,y\in\Mani\\\dist_\Mani(x,y)\geq b_n}}  \abs{\frac{\widehat C_{n,\infty}(x,y)}{C_\infty(x,y)} -1} \ge \eps \right) \leq e^{-\const_0n^{\const_1}}. 
    \end{equation}
\end{theorem}

\begin{proof}%[Proof of Theorem \ref{th:hatCnInfVSCInf}]
As this proof is slightly longer, we will structure it in paragraphs. 
    We start proving that $C_\infty(x,y)$ is a high-probability lower-bound for $\widehat C_{n,\infty}(x,y)$. 
    
    \vspace{0,1cm}
    
    \noindent{\bf High-probability lower-bound.} For any $x,y\in\Mani$ let $\hatgeod(x,y)=(u_i)_{i=1}^\ell$, $\ell=\widehat K_{n,\infty}(x,y)$, be a $\widehat C_{n,\infty}(x,y)$ geodesic. Also let $\cgeodi{i}$ be 
    the geodesic attaining $D_{\infty,\alpha}(u_i,u_{i+1})$
    and consider the concatenation
    $\cpath_\hatgeod = \cgeodi{1} + \cdots + \cgeodi{\ell-1}$.
    %be its $D_{\infty,\alpha}$-continuous extension (see Definition \ref{def:discretPathContinuous}).
    Let also $\cgeodi{\ell}$ be the $D_{\infty,\alpha}$-geodesic from $u_\ell$ to $y$, where we recall
    that by definition of $\mathcal U_n,$ $\dist_\Mani(u_\ell,y)\leq \delta_n.$
    We compute
    \begin{align}\label{eq:uppberBoundCinfty1}
    C_\infty(x,y) &\leq
    \int_{\cpath_\hatgeod} e^{-\frac{A(\length_{\cpath_\hatgeod}(z))}{d}}f(z)^{-\frac{\alpha-1}{d}} d\cpath_\hatgeod
    + \int_{\cgeodi{\ell}} e^{-\frac{A\left(\length\cpath_\hatgeod + \length_{\cgeodi{\ell}}(z)\right)}{d}}f(z)^{-\frac{\alpha-1}{d}} d\cgeodi{\ell}\nonumber\\
    &\leq \sum_{i=1}^{\ell-1}  \int_{\cgeodi{i}} e^{-\frac{A\left(\sum_{j=1}^{i-1} \length\cgeodi{j} + \length_{\cgeodi{i}}(z)\right)}{d}}f(z)^{-\frac{\alpha-1}{d}} d\cgeodi{i}
    + (\minv f)^{-\frac{\alpha-1}{d}} \length \cgeodi{\ell}.
    \end{align}
    Using 
    Proposition \ref{prop:boundedCurvatureDinftyGeod} and
    \eqref{eq:geodLengthSmallApprox} we have the bound $\length \cgeodi{\ell} \leq \dist_\Mani(u,y)(1+\lo(1))\leq \delta_n(1+\lo(1))$, where the error term $\lo(1)$ holds uniformly on $x,y$ as $\delta_n\to0$. Thus for the second term in the r.h.s. of \eqref{eq:uppberBoundCinfty1} we obtain
    \begin{align*}
    (\minv f)^{-\frac{\alpha-1}{d}} \length \cgeodi{\ell}& \leq (\minv f)^{-\frac{\alpha-1}{d}}
    \delta_n(1+\lo(1)).
    \end{align*}
    On the other hand, note that by the definition of $\mathcal U_n(x,y)$ (in \eqref{eq:defUn}) and the fact that $\length\cgeodi{j}\geq \dist_\Mani(u_j,u_{j+1})\geq \delta_n$, we have
    $
    \sum_{j=1}^{i-1} \length\cgeodi{j} \geq  (i-1)\delta_n
    $
    so that
    $$
    \forall 1\leq i \leq \ell, z\in \cgeodi{i},\quad A\left(\sum_{j=1}^{i-1} \length\cgeodi{j} + \length_{\cgeodi{i}}(z)\right) \geq A(i\delta_n) -\modulus_A(\delta_n).
    $$
    This implies, for the first term in the r.h.s. of \eqref{eq:uppberBoundCinfty1},
    \begin{align*}
    &\sum_{i=1}^{\ell-1}  \int_{\cgeodi{i}} e^{-\frac{A\left(\sum_{j=1}^{i-1} \length\cgeodi{j} + \length_{\cgeodi{i}}(z)\right)}{d}}f(z)^{-\frac{\alpha-1}{d}} d\cgeodi{i}\\
    &\leq e^{\frac{\modulus_A(\delta_n)}{d}} \sum_{i=1}^{\ell-1} e^{-\frac{A\left(i\delta_n\right)}{d}} \int_{\cgeodi{i}} f(z)^{-\frac{\alpha-1}{d}} d\cgeodi{i}
    =e^{\frac{\modulus_A(\delta_n)}{d}} \sum_{i=1}^{\ell-1} e^{-\frac{A\left(i\delta_n\right)}{d}} D_{\infty,\alpha}(u_i,u_{i+1})\\
    &=e^{\frac{\modulus_A(\delta_n)}{d}} \widehat C_{n,\infty}(x,y).
    \end{align*}
    Combining the two bounds we obtain
    $$ 
    C_\infty(x,y) \leq e^{\frac{\modulus_A(\delta_n)}{d}} \widehat C_{n,\infty}(x,y) + (\minv f)^{-\frac{\alpha-1}{d}}
    \delta_n(1+\lo(1)).
    $$
    Since $C_\infty(x,y)\geq e^{-\frac{\maxv A}{d}}(\maxv f)^{-\frac{\alpha-1}{d}}b_n$ whenever $\dist_\Mani(x,y)\geq b_n$,
    we obtain 
    \begin{equation}\label{eq:limitIntLower}
    {\inf_{\substack{x,y\in\Mani\\\dist_\Mani(x,y)\geq b_n}}} \frac{\widehat C_{n,\infty}(x,y)}{C_\infty(x,y)} \asgeq 
    e^{-\frac{\modulus_A(\delta_n)}{d}} -  \frac{e^{\frac{\maxv A}{d}}(\minv f)^{-\frac{\alpha-1}{d}} (1+\lo(1))}{(\maxv f)^{-\frac{\alpha-1}{d}}}\frac{\delta_n}{b_n}.
    \end{equation}
    The r.h.s. converges  to $1$ uniformly on $x,y$ as $n\to\infty.$\par
This concludes the first step of our proof.     
    We will now prove that $C_\infty(x,y)$ is a high-probability upper-bound for $\widehat C_{n,\infty}(x,y)$.  Recall the set $\mathcal U_n(x,y)$ in \eqref{eq:defUn}. The quantity $\widehat C_{n,\infty}(x,y)$ being defined as
    an infimum over $\mathcal U_n(x,y)$, we need only prove that, with high probability and for fixed $\eps>0$ and any $x,y\in\Mani$ such that $\dist_\Mani(x,y)\geq b_n$, there exists a mesoscopic sample path $\path(x,y)=(u_1,\cdots,u_\ell)\in \mathcal U_n(x,y)$ of size $1\leq \ell(x,y)\leq n$ such that
    $$
    \frac{\sum_{i=1}^{\ell(x,y)-1} e^{-\frac{A\left(i\delta_n\right) }{d}}   
     D_{\infty,\alpha}(u_i,u_{i+1})}{C_\infty(x,y)} \leq 1 + \eps.
    $$\par

\vspace{0,1cm}
    
    \noindent{\bf Construction of a mesoscopic path.}
    Let us construct such a path $\path\equiv \path(x,y)$. 
    Let $\eps>0$ be arbitrary and consider a curve $\cgeod\equiv \cpath^*_{x,y}$  satisfying
    \begin{equation}\label{eq:curveApproxCinfty}
    \frac{\int_0^{\length\cgeod} e^{-\frac{A(t)}{d}} f(\cgeod(t))^{-\frac{\alpha-1}{d}}dt}{C_\infty(x,y)} = 1+\eps,
    \end{equation}
    where for simplicity we have assumed that $\cgeod$ is parametrized with unit-speed. Furthermore, by Proposition \ref{prop:restrictionToBoundedCurvature}, we may assume that $\cgeod\in\mathcal F_{\kappa^*}$ and $\length\cgeod\leq \const_\length\dist(x,y)\leq \const_\length\diameter(\Mani)$.\par
    Define a  partition 
    $\partition=(\cgeod(t_i))_{i=1}^{\ell}$ of $\cgeod$ similarly as in the proof of Proposition \ref{prop:restrictionToBoundedCurvature} 
    but replacing $\delta$ by $\delta_n+2\Delta_n$.
    Note that $\cgeod$ being continuous, we have
    $$
    \forall 1\leq i\leq \ell-1\colon\quad \dist_\Mani(\cgeod(t_i),\cgeod(t_{i+1}))=\delta_n+2\Delta_n.
    $$
    Also, since  
    $$\forall 1\leq i\leq \ell-1\colon\quad \length_\cgeod(t_{i+1})-\length_\cgeod(t_{i})\geq \dist_\Mani(\cgeod(t_i),\cgeod(t_{i+1}))=\delta_n+2\Delta_n,$$
    we have
    \begin{equation}\label{eq:boundEll_uppberboundHatCninfty}
    (\ell-1) (\delta_n+2\Delta_n)\leq \length\cgeod  \leq \const_\length\diameter(\Mani).
    \end{equation}
    In particular, by Assumption \ref{cond:DeltanBound}, $(\ell-1) \leq \const_\length\diameter(\Mani)  n^{\frac{1-\eta}{d}}$ so that
    \begin{equation}\label{eq:boundEll_belongUn}
    \ell \leq n
    \end{equation}
    for large enough $n$.
    Furthermore, 
    \begin{equation}\label{eq:mesPathConstruction1}
    \forall 1\leq i\leq \ell-1, \quad\forall \substack{z\in \ball(\cgeod(t_{i}),\Delta_n )\\z'\in \ball(\cgeod(t_{i+1}),\Delta_n )}
    \colon \quad \delta_n\leq \dist_\Mani(z,z')\leq \delta_n + 4\Delta_n.
    \end{equation}
    Indeed, this is a simple consequence of the reverse triangle inequality applied twice, which yields, for such pairs $z,z'$:
    \begin{align*}
    \dist_\Mani(z,z') &= \dist_\Mani(z, \cgeod(t_{i+1})) \pm \Delta_n = \dist_\Mani(\cgeod(t_{i}), \cgeod(t_{i+1})) \pm 2\Delta_n
    \\&
    =\delta_n+2\Delta_n\pm2\Delta_n.
    \end{align*}
    Let us then construct the mesoscopic path $\path\equiv \path(x,y)=(u_1,\cdots,u_\ell)$ by letting $u_1=x$ and, for $2\leq i\leq \ell(x,y)$, letting $u_i$ to be any point in $Q_n \cap \bar\ball(\cgeod({t_{i}}),\Delta_n )$ (by
    the definition of $\Delta_n$ this always exists). Finally, if $\dist_\Mani(u_\ell,y)\geq \delta_n$ we consider 
    $\path=(u_1,\cdots,u_\ell,y)$ instead.  
    Note that, from \eqref{eq:boundEll_belongUn} and also \eqref{eq:mesPathConstruction1} 
    %together with $\theta<1$, 
    we have $\path(x,y)\in \mathcal U_n(x,y).$
    \par

    With our path in hands, we can pass to the last step of our proof which is to prove  that $C_\infty(x,y)$ is a high-probability upper-bound for $\widehat C_{n,\infty}(x,y)$.

   \vspace{0,1cm}
    
    \noindent{\bf High-probability upper-bound.}
    Let $t_{\ell+1}=\length\cgeod$ so that $\cgeod(t_{\ell+1})=y$. Having defined the partition $(\cgeod(t_{i}))_{i_1}^{\ell+1}$ 
    and the discrete path $\path(x,y)$,  the idea is to prove that the following inequalities \eqref{eq:CinftyLowBound1}
    and \eqref{eq:CinftyLowBound2} hold with high probability
    for arbitrary $\eps_1,\eps_2>0$, and
    uniformly on the set $\{x,y\in\Mani\colon \dist_\Mani(x,y)\geq \delta_n\}$ and on the choice of $\cgeod(x,y)$. Given \eqref{eq:curveApproxCinfty}, it is enough to prove that
    \begin{align}
    \int_0^{\length\cgeod} e^{-\frac{A(t)}{d}} f(\cgeod(t))^{-\frac{\alpha-1}{d}}dt &\geq e^{-\modulus_A(\eps_1)} \sum_{i=1}^{\ell(x,y)-1}e^{-\frac{A(i\delta_n)}{d}}D_{\infty,\alpha}(\cgeod(t_{i}),\cgeod(t_{i+1})) 
    \label{eq:CinftyLowBound1}
    \end{align}
    and also that
    \begin{align}
    \label{eq:CinftyLowBound2}
    &\sum_{i=1}^{\ell(x,y)-1}e^{-\frac{A(i\delta_n)}{d}}D_{\infty,\alpha}(\cgeod(t_{i}),\cgeod(t_{i+1})) 
    \nonumber\\
    &\geq(1-\eps_2)\sum_{i=1}^{\ell(x,y)-1} e^{-\frac{A\left(i\delta_n\right) }{d}}   
     D_{\infty,\alpha}(u_i,u_{i+1}).
    \end{align}
    The first inequality holding in a high-probability event $\mathcal A_{n,\eps_1}$, and the second in another high-probability event $\mathcal B_{n,\eps_2}$.
    Here we recall that, by Assumption \eqref{cond:functionA}, the function $A$ is uniformly continuous so that 
    for arbitrary $\eps_1'>0$ and choosing $\eps_1$ small enough we have $e^{-\modulus_A(\eps_1)}\geq(1-\eps_1')$.\par
    %let us first compare $C_\infty(x,y)$ against the following sum: $$\sum_{i=1}^{\ell-1}e^{-\frac{A(k\delta_n)}{d}}D_{\infty,\alpha}(\cgeod(t_{i}),\cgeod(t_{i+1})).$$
    We first prove that \eqref{eq:CinftyLowBound1} holds in a high probability event $\mathcal A_{n,\eps_1}$. 
    Letting  $\cpath_i\coloneqq \restr{\cpath_{x,y}^*}{[t_i,t_{i+1}]}$, and 
    following the proof of Proposition \ref{prop:restrictionToBoundedCurvature}
    with $\kappa\equiv\kappa(x,y)$ therein replaced by the universal constant $\kappa^*$, we obtain
    $$
    \sup_{\substack{x,y\in\Mani\\\dist_\Mani(x,y)\geq \delta_n}} \sup_{1\leq i\leq \ell(x,y)-1}
     \abs{ \frac{\sum_{j=1}^{i} \length\cpath_j}{i(\delta_n+2\Delta_n) } -1 }\leq \eps_1'
    $$
    on the event $\{\delta_n+2\Delta_n\leq \delta'\}$
    where  $\delta'>0$ is taken small enough to ensure that then $\bO$-term in the r.h.s. of \eqref{eq:geodLengthSmallApprox} is smaller than $\eps_1'=\eps_1/(2\const_\length\diameter(\Mani))$.
    From the latter, 
    \begin{align}\label{eq:limitIntUpper_Aapprox}
    \sup_{\substack{x,y\in\Mani\\\dist_\Mani(x,y)\geq \delta_n}} 
    \sup_{1\leq i\leq \ell(x,y)-1} A\left(\sum_{j=1}^{i} \length\cpath_j \right) 
    &\leq A\Big(i(\delta_n+2\Delta_n) (1+ \eps_1')\Big)
    \nonumber\\
    &\leq A\Big(i\delta_n\Big) 
    + \modulus_A\Big(\ell2\Delta_n(1+\eps_1')+\ell\delta_n\eps_1'\Big).
    %\nonumber\\
    %& \leq A\left(i\delta_n\right) 
    %+ \modulus_A\left(\const_\length\diameter(\Mani)\frac{\Delta_n}{\delta_n} + \const_\length\diameter(\Mani) \lo(1)\right).
    \end{align}
    Using \eqref{eq:boundEll_uppberboundHatCninfty} we have, for the argument of the modulus $\modulus_A$ in the last line above,
    $$
    \ell2\Delta_n(1+\eps_1')+\ell\delta_n\eps_1'\leq \frac{2\Delta_n}{\delta_n}\left(\const_\length\diameter(\Mani)+\eps_1\right)+ \frac{\eps_1}{2}+\Delta_n(1+\eps_1')+\delta_n\eps_1'.
    $$
    Plugging this bound in \eqref{eq:limitIntUpper_Aapprox} and the latter in \eqref{eq:CinftyLowerBound1_1} we obtain
    that \eqref{eq:CinftyLowBound1} holds uniformly on $x,y\in\Mani$ in the event $\mathcal A_{n,\eps_1}\coloneqq\{\delta_n+\Delta_n\leq\delta'\}\cap\left\{ \frac{2\Delta_n}{\delta_n}\left(\const_\length\diameter(\Mani)+\eps_1\right)+ \Delta_n(1+\eps_1')+\delta_n\eps_1'\leq \eps_1/2\right\}$. \par
    \par
    Let us now prove that \eqref{eq:CinftyLowBound2} holds in a high-probability event $\mathcal B_{n,\eps_2}$.
    Note that by the reverse triangle inequality applied twice we have
     \begin{align*}
    1\leq i\leq \ell-1,\quad 
    D_{\infty,\alpha}(\cgeod(t_{i}),\cgeod(t_{i+1})) = D_{\infty,\alpha}(u_i,u_{i+1}) \pm 2\Delta_n.
     \end{align*}
     Thus
         \begin{align*}
    \forall 1\leq i\leq \ell-1\colon\quad 
    &\abs{\frac{D_{\infty,\alpha}(\cgeod(t_{i}),\cgeod(t_{i+1}))}{D_{\infty,\alpha}(u_i,u_{i+1})} - 1}
    \leq \frac{2\Delta_n}{D_{\infty,\alpha}(u_i,u_{i+1})}.
     \end{align*}\par
    Using Lemma \ref{le:metricDerivativeOfCurves}
    and a first order Taylor expansion of the integral
    we obtain, for any unit-speed curve $\cpath\in\mathcal{F}_\kappa$,
\begin{align*}
     \frac{\int_s^t  f(\cpath(u))^{-\frac{\alpha-1}{d}} du}{f(\cpath(s))^{-\frac{\alpha-1}{d}}\dist_\Mani(\cpath(t),\cpath(s))} 
     &=\frac{ f(\cgeod(s))^{-\frac{\alpha-1}{d}}(t-s) + \bO(\abs{t-s}^2)}{f(\cpath(s))^{-\frac{\alpha-1}{d}}\dist_\Mani(\cpath(t),\cpath(s))} 
     \\&= 1 + \bO(\abs{t-s}),
\end{align*}
where the error term $\bO(\abs{t-s})$ is uniform over the entire class $\mathcal{F}_\kappa$.\par
    In particular, as $n\to\infty$ (so that
    $\delta_n\to0$), we have
    $D_{\infty,\alpha}(u_i,u_{i+1})=f(u_i)^{-\frac{\alpha-1}{d}}\dist_\Mani(u_i,u_{i+1})(1+\lo(1))$, where the error term $\lo(1)$ is uniform
    on $1\leq i\leq \ell(x,y)-1$, and also on $x,y\in\Mani$ and $\cgeod(x,y)$.  Since by construction $\dist_\Mani(u_i,u_{i+1})\geq\delta_n$ we have,
    for $n$ large enough to ensure $\lo(1)<1/2$, that
    \begin{equation*}
    \sup_{\substack{x,y\in\Mani\\\dist_\Mani(x,y)\geq \delta_n}} \sup_{1\leq i\leq \ell(x,y)-1}
    \abs{ \frac{D_{\infty,\alpha}(\cgeod(t_{i}),\cgeod(t_{i+1}))}{D_{\infty,\alpha}(u_i,u_{i+1})} - 1}\leq \frac{2}{\maxv f^{-\frac{\alpha-1}{d}}(1-1/2)}\frac{\Delta_n}{\delta_n}.
    \end{equation*}
    Thus, on the event $\mathcal{B}_{n,\eps_2}\coloneqq\left\{\frac{4}{\maxv f^{-\frac{\alpha-1}{d}}}\frac{\Delta_n}{\delta_n}\leq\eps_2\right\}$, we have
    \begin{equation*}
    \forall \substack{x,y\in\Mani\\\dist_\Mani(x,y)\geq \delta_n}\colon\quad
   \frac{\sum_{i=1}^{\ell(x,y)-1} e^{-\frac{A\left(i\delta_n\right) }{d}}   
     D_{\infty,\alpha}(\cpath^*_{x,y}(t_{i}),\cpath^*_{x,y}(t_{i+1}))}{\sum_{i=1}^{\ell(x,y)-1} e^{-\frac{A\left(i\delta_n\right) }{d}}   
     D_{\infty,\alpha}(u_i,u_{i+1}))} = 1\pm \eps_2.
    \end{equation*}
    In particular
    \begin{equation}\label{eq:limitIntUpper2}
        \inf_{\substack{x,y\in\Mani\\\dist_\Mani(x,y)\geq \delta_n}} 
        \frac{\sum_{i=1}^{\ell(x,y)-1} e^{-\frac{A\left(i\delta_n\right) }{d}}   
     D_{\infty,\alpha}(\cpath^*_{x,y}(t_{i}),\cpath^*_{x,y}(t_{i+1}))}{\sum_{i=1}^{\ell(x,y)-1} e^{-\frac{A\left(i\delta_n\right) }{d}}   
     D_{\infty,\alpha}(u_i,u_{i+1}))}\geq 1 - \eps_2.  
    \end{equation}
    This proves \eqref{eq:CinftyLowBound2} on the high probability event $\mathcal B_{n,\eps_2}$.\par
    Combining \eqref{eq:CinftyLowBound1} and \eqref{eq:CinftyLowBound2} we obtain, 
    on the event $\mathcal A_{n,\eps_1} \cap \mathcal B_{n,\eps_2}$, that
    \begin{align*}
     &\sup_{\substack{x,y\in\Mani\\\dist_\Mani(x,y)\geq \delta_n}} \frac{\widehat C_{n,\infty}(x,y)}{C_\infty(x,y)}
    \leq\sup_{\substack{x,y\in\Mani\\\dist_\Mani(x,y)\geq \delta_n}} \frac{\sum_{i=1}^{\ell(x,y)-1} e^{-\frac{A\left(i\delta_n\right) }{d}}   
     D_{\infty,\alpha}(u_i,u_{i+1}))}{C_\infty(x,y)} 
    \end{align*}     
    is bounded by $e^{{\modulus_A(\eps_1)}}(1-\eps_2)^{-1}$.
    \par
    
       \vspace{0,1cm}
    
    \noindent{\bf Conclusion of the proof of Theorem \ref{th:hatCnInfVSCInf}.}
    Combining the above with \eqref{eq:limitIntLower} we obtain, for arbitrary but fixed ${\eps}>0$ and for large enough $n\geq1$ so that
    the r.h.s. of \eqref{eq:limitIntLower} is larger than $1-\eps$, 
    \begin{align*}
    &\P\left(\forall  \substack{x,y\in\Mani\\\dist_\Mani(x,y)\geq b_n}\colon \quad  1-\eps\leq \frac{\widehat C_{n,\infty}(x,y)}{C_\infty(x,y)}  \leq e^{{\modulus_A(\eps)}}(1-\eps)^{-1} \right)\\
    &\geq \P\left(\mathcal A_{n,\eps} \cap \mathcal B_{n,\eps}\right). 
    \end{align*}
    By Assumption \ref{cond:DeltanBound}
    and Lemma \ref{le:DeltanBound} the latter probability is lower-bounded by $1-e^{-\const_0n^{\const_1}}$ for large enough $n$ and some $\const_0,\const_1$ depending on 
    ${\eps}$ (and also on $\maxv f^{-\frac{\alpha-1}{d}}$ and $\const_\length\diameter(\Mani)$). Recalling that, by Assumption \ref{cond:functionA}, $\modulus_A(\eps)\to0$ as $\eps\to0$, the above yields \eqref{eq:hatCnInftyVSCinfty}.
\end{proof}

\section{Final step: Proofs of the main results}\label{sec:final_step}
% \subsection{Proofs of manifold FPP results}\label{sec:ProofFermatCostConvergence}

In this section we collect the results of sections \ref{sec:macroscopicLimit} and \ref{sec:mesoscopicApproximation} to prove Theorems \ref{th:FermatCostConvergence} and  \ref{th:GHconvergence2}. {These proofs will be done for $C_n=C_n^{\dist_\Mani}$. We remind the reader that, according to Proposition \ref{prop:equivalentModels}, the same statements will hold  for $C_n^{\dist_E}$}.
\begin{proof}[Proof of Theorem \ref{th:FermatCostConvergence}]
The Theorem easily follows from writing 
    \begin{align*}
\abs{\frac{n^{\frac{\alpha-1}{d}} C_n(x,y)}
     {\mu(d,\alpha)C_{\infty}(x,y) } -1 } &\leq 
     \abs{\frac{n^{\frac{\alpha-1}{d}}C_n(x,y)}{{\mu(d,\alpha)}\widehat C_{n,\infty}(x,y)} - 1}\abs{\frac{{\mu(d,\alpha)}\widehat C_{n,\infty}(x,y)}{\mu(d,\alpha)C_{\infty}(x,y)} -1 } 
     \\&+  \abs{\frac{n^{\frac{\alpha-1}{d}}C_n(x,y)}{{\mu(d,\alpha)}\widehat C_{n,\infty}(x,y)} - 1} +\abs{\frac{{\mu(d,\alpha)}\widehat C_{n,\infty}(x,y)}{\mu(d,\alpha) C_{\infty}(x,y)} -1 }, 
    \end{align*}   
    and then using Theorems \ref{theo:hatCnVshatCninfty} and \ref{th:hatCnInfVSCInf} to bound each term with high probability.
\end{proof}

The following Lemma, which is key to Theorem \ref{th:GHconvergence2}, has an interest by itself, as it states that the random  costs $n^{\frac{\alpha-1}{d}}C_n(x,y)$ converge uniformly on all $x,y\in\Mani$ towards  $\mu(d,\alpha)C_\infty(x,y)$. 

\begin{lemma}\label{le:GHconvergence1}
 Let $\eps>0$. There exists constants $\const_0>0,\const_1>0$ 
such that for large enough $n$,
\begin{equation}
	\P\left(\sup_{\substack{x,y\in\Mani}}  \abs{n^{\frac{\alpha-1}{d}}C_{n}(x,y) -\mu(d,\alpha)C_\infty(x,y)} \ge \eps \right) \leq e^{-\const_0n^{\const_1}}. 
\end{equation}
\end{lemma}
The convergence in Theorem \ref{th:FermatCostConvergence} gives a tighter control on the error 
$\abs{n^{\frac{\alpha-1}{d}}C_{n}(x,y) -\mu(d,\alpha)C_\infty(x,y)}$ since it scales
with $C_\infty(x,y)$, whereas in the above Lemma \ref{le:GHconvergence1} it does not. We are able to obtain the uniform result in the latter by loosening the control on this error, specially when we consider
points $x_n,y_n\in\Mani$ such that $\dist_\Mani(x_n,y_n)\to 0$ for which we also have $C_\infty(x_n,y_n)\to 0$.
\begin{proof}[Proof of Lemma  \ref{le:GHconvergence1}]
    Take $b_n\toas{n\to\infty} 0$ satisfying Assumption \ref{cond:bn}. Note that, from Assumption \ref{cond:DeltanBound}, in particular we have
    $b_n\geq Cn^{-\frac{1-\eta\alpha}{d}}$ where we recall that $\eta\alpha<1$ (Assumption \ref{cond:DeltanBound}).
	Given Theorem \ref{th:FermatCostConvergence} it suffices to prove
	\begin{equation}\label{eq:unifConvSmallDist}
    \P\left( \sup_{\substack{x,y\in\Mani\\\dist_\Mani(x,y)\leq b_n}} n^{\frac{\alpha-1}{d}}C_n(x,y) + \mu(d,\alpha)C_\infty(x,y) >\eps  \right)\leq e^{-\const_0n^{\const_1}}
	\end{equation}
	for large enough $n$. On the one hand we clearly have, with probability $1$, that
	$
	\sup_{\substack{x,y\in\Mani\\\dist_\Mani(x,y)\leq b_n}} C_\infty(x,y) \leq \minv f^{-\frac{\alpha-1}{d}} b_n
    \toas{n\to\infty}0.
	$
    On the other hand, by Lemmas \ref{le:JnAsmp} and \ref{le:DeltanBound}, we have $C_n(x,y)\leq D_{n,\alpha}(x,y)$ with probability at least $1-e^{-\const_0n^{\const_1}}$ for some positive constants $\const_0,\const_1$. 
    Take $q\in Q_n$ satisfying $b_n\leq \dist_\Mani(x,q)\leq 3b_n$ and  $b_n\leq \dist_\Mani(q,y)\leq 3b_n$. By Lemma
    \ref{le:DeltanBound} this can be done with probability at least $1-e^{-\const_0'n^{\const_1'}}$. Then, since $q\in Q_n$, using 
    the triangle inequality for $D_{n,\alpha}$, and then \eqref{eq:damelin_concentration} for each resulting term, we obtain with probability at least $1-e^{-\const_0''n^{\const_1''}}$
    that, 
    $n^{\frac{\alpha-1}{d}}C_n(x,y)\leq n^{\frac{\alpha-1}{d}}D_{n,\alpha}(x,y)\leq 2(\mu(d,\alpha)+1)(\minv f)^{-\frac{\alpha-1}{d}}3b_n\toas{n\to\infty}0$.
    Combining both bounds we obtain \eqref{eq:unifConvSmallDist}.
\end{proof}

\begin{proof}[Proof of Theorem \ref{th:GHconvergence2}]
	We first note that $C_\infty$ satisfies the triangle inequality, this is an easy consequence of the fact that the function $e^{-\frac{A(s)}{d}}$ is
	non-increasing.
	We now adapt the argument leading to Theorem 7 in \cite{FernandezBorghiniMindlinGroisman2023} to prove the  theorem.
	Let us write 
    $\tilde C_n = n^{\frac{\alpha-1}{d}}C_n$ and $\tilde C_\infty =\mu(d,\alpha) C_\infty.$ The distance $\dist_{\mathcal N}$, an infimum over correspondence sets, is upper-bounded 
    by picking the particular correspondence set $A\subset Q_n \times\Mani$ given by
    $$A=\{ (q,q) \colon  q \in Q_n\} \cup \{(q_y,y) \colon y\in\Mani, \dist_\Mani(q_y,y)=\dist_\Mani(Q_n,y)\}.$$
	For this correspondence set we compute
	\begin{align*}
	&2\dist_{\mathcal{N}}((Q_n,C_n),(\Mani,\tilde C_\infty))
	\\&\leq  \sup_{q,q'\in Q_n} \abs{\tilde C_n(q,q')-\tilde  C_\infty(q,q')} \vee  \sup_{y,y'\in \Mani} 
	 \abs{\tilde C_n(q_y,q_{y'})-\tilde  C_\infty(y,y')} \\
	 &\leq 2\sup_{x,y\in \Mani} \abs{\tilde C_n(x,y)-\tilde C_\infty(x,y)} \vee  \sup_{y,y'\in \Mani} 
	  \abs{\tilde C_\infty(q_y,q_{y'})-\tilde C_\infty(y,y')}.
	\end{align*}
	The first supremum in the r.h.s. is bounded with high probability using Lemma \ref{le:GHconvergence1}. 
    For the second, note that, using the triangle inequality of $C_\infty$ (twice), for any $x,y,z\in\Mani$, we have
    $
    \abs{\tilde C_\infty(x,y) -  \tilde C_\infty(z,y)} \leq \max\{\tilde C_\infty(x,z) ,\tilde C_\infty(z,x) \}.
    $
    Thus
    \begin{align*}
    \abs{\tilde C_\infty(q_y,q_{y'})-\tilde C_\infty(y,y')} &\leq \max\{\tilde C_\infty(y,q_y) ,\tilde C_\infty(q_y,y) \} \\&+ \max\{\tilde C_\infty(y',q_{y'}),\tilde C_\infty(q_{y'},y')\}.
    \end{align*}
    Since, for any pair $y,q_y$ we have
    $$
       \max\{\tilde C_\infty(y,q_y), \tilde C_\infty(q_y,y)\} \leq \mu(d,\alpha)\minv f^{-\frac{\alpha-1}{d}}\Delta_n,
    $$
    taking supremum we obtain
    $$
    \sup_{y,y'\in \Mani} 
    \abs{\tilde C_\infty(q_y,q_{y'})-\tilde C_\infty(y,y')}  \leq 2\mu(d,\alpha)\minv f^{-\frac{\alpha-1}{d}}\Delta_n.
    $$
    The latter is bounded with high probability using Lemma \ref{le:DeltanBound}. Combining both bounds we conclude the proof.
\end{proof}

\appendix

\section[Step 0: $C^{\dist_\Mani}_n\approx C^{\dist_E}_n$]{Step 0: Asymptotic equivalence between $C^{\dist_\Mani}_n$ and $C^{\dist_E}_n$}\label{sec:step0}
We begin with the following result concerning the geometry of $\Mani\subset \R^D$, which compares the Euclidean and the
geodesic distance between points that are close together.
\begin{proposition}[Corollary 4 \cite{Bernstein2000}, Lemma 3 \cite{Boissonnatetal2019}]\label{prop:euclideaVSgeodesic}
	Let $\Mani$ be a smooth compact manifold embedded in $\R^D$. Define the functions
	$$
	\zeta_0(\eps) \coloneqq \sup_{\substack{x,y\in\Mani\subset \R^D\\ \dist_E(x,y)\leq \eps}} \abs{\frac{\dist_\Mani(x,y)}{\dist_E(x,y)} - 1}    ,\text{ and}\quad \zeta(\eps) \coloneqq \sup_{\substack{x,y\in\Mani\subset \R^D\\ \dist_\Mani(x,y)\leq \eps}} \abs{\frac{\dist_\Mani(x,y)}{\dist_E(x,y)} - 1}.    
	$$
Then we have $0\leq \zeta_0(\eps)  \leq \zeta(\eps) = \lo(1)$ where error term $\lo(1)$ holds as $\eps\to0$. 
    \end{proposition}

The asymptotic equivalence between $C^{\dist_\Mani}_n$ and $C^{\dist_E}_n$ is a consequence of Proposition \ref{prop:euclideaVSgeodesic} and the fact that
$C^{\dist_\Mani}_n$ and $C^{\dist_E}_n$-geodesic paths only have small jumps, as stated in Lemma \ref{le:JnAsmp}.

\begin{proof}[Proof of Proposition \ref{prop:equivalentModels}]
	For this proof, we will use Theorem \ref{th:FermatCostConvergence} with $\dist=\dist_\Mani$, which is itself proved using arguments that are independent from this section. Consider $C^{\dist_\Mani}_{n,1}$ and  $C^{\dist_\Mani}_{n,2}$ defined as $C^{\dist_\Mani}_n$ but using respectively  $\delta_n=\delta_{n,1}$ and $\delta_n=\delta_{n,2}$. If both $\delta_{n,1}$ and $\delta_{n,2}$ satisfy
	Assumption \ref{cond:DeltanBound} then, by Theorem \ref{th:FermatCostConvergence} with $\dist=\dist_\Mani$, for any $\eps>0$ there exists $\const_0,\const_1>0$ such that,
	for sufficiently large $n$,
    \begin{align*}
	&\P\left(\sup_{\substack{x,y\in\Mani\\\dist_\Mani(x,y)\geq b_n}} \abs{ \frac{n^{\frac{\alpha-1}{d}}C^{\dist_\Mani}_{n,1}(x,y)}{\mu(d,\alpha)C_{\infty}(x,y)} - 1}
	        \vee \sup_{\substack{x,y\in\Mani\\\dist_\Mani(x,y)\geq b_n}} \abs{ \frac{\mu(d,\alpha)C_{\infty}(x,y)}{n^{\frac{\alpha-1}{d}}C^{\dist_\Mani}_{n,2}(x,y)} - 1}
	       \geq \eps\right)
	\\& \leq 2e^{-\const_0n^{\const_1}}.
	\end{align*}
    Thus, for sufficiently large $n$ and possibly different constants $\const_0>0,\const_1>0$ we have
	\begin{equation}\label{eq:equivalenceC_n,delta_n}
	\P\left( \sup_{\substack{x,y\in\Mani\\\dist_\Mani(x,y)\geq b_n}} \abs{ \frac{C^{\dist_\Mani}_{n,1}(x,y)}{C^{\dist_\Mani}_{n,2}(x,y)} - 1}\geq \eps \right)
	\leq e^{-\const_0n^{\const_1}}.
	\end{equation}
    \par
    The idea now will be to compare $C^{\dist_E}_n$ with the costs $C^{\dist_\Mani}_n,C^{\dist_\Mani}_{n,1},$ and $C^{\dist_\Mani}_{n,2}$ which are computed using, respectively, $\delta_n=\delta_n^{\dist_E}$, $\delta_{n,1}=2\delta_n^{\dist_E}$, and $\delta_{n,2}=\delta_n^{\dist_E}$. Specifically, we will prove that 
    \begin{align}\label{eq:equivalence1}
		&
        \sup_{x,y\in\Mani} \frac{C^{\dist_E}_n(x,y)}{C^{\dist_\Mani}_{n,1}(x,y)}\asleq 1
		%\nonumber\\&\leq\P\left( \left(\frac{1}{1+\zeta(\Delta_n/\theta)}\right)^{\alpha+\maxv A}; \frac{\Delta_n}{\theta}<1 \right) 
		%\leq e^{-\const_0 n^{\const_1}}
	\end{align}
    and that for any $\epsilon>0$ and sufficiently large $n$
    \begin{equation}\label{eq:equivalence2}
		\P\left({\inf_{x,y\in\Mani}} \frac{C^{\dist_E}_n(x,y)}{C^{\dist_\Mani}_{n,2}(x,y)}{<} 1-\eps\right)\leq e^{-\const_0 n^{\const_1}}.
	\end{equation}
	Having proved the above, the proposition follows 
    by writing
    \begin{align*}
    \frac{C^{\dist_E}(x,y)}{C^{\dist_\Mani}_n(x,y)}
    \eq \frac{C^{\dist_E}(x,y)}{C^{\dist_\Mani}_{n,1}(x,y)} \frac{C^{\dist_\Mani}_{n,1}(x,y)}{C^{\dist_\Mani}_{n}(x,y)}
    \text{ and }
    \frac{C^{\dist_E}(x,y)}{C^{\dist_\Mani}_n(x,y)}\eq \frac{C^{\dist_E}(x,y)}{C^{\dist_\Mani}_{n,2}(x,y)} \frac{C^{\dist_\Mani}_{n,2}(x,y)}{C^{\dist_\Mani}_{n}(x,y)}
    \end{align*}
    and then combining \eqref{eq:equivalence1} and \eqref{eq:equivalence2}, together with   \eqref{eq:equivalenceC_n,delta_n}, to bound each term with high probability. 
So, to conclude the proof of the proposition, it only remains to prove \eqref{eq:equivalence1} and \eqref{eq:equivalence2}. 
Let us start proving \eqref{eq:equivalence1}. Note that by Proposition \ref{prop:euclideaVSgeodesic} and for $n\geq1$ large enough to ensure $\zeta_0(\delta_n^{\dist_E})<1$ we have
	$$
	\forall q,q'\in\Mani\colon\quad \dist_E(q,q')< \delta_n^{\dist_E} \implies \dist_\Mani(q,q')\leq \delta_n^{\dist_E}(1+\zeta_0(\delta_n^{\dist_E})) <2\delta_n^{\dist_E}.
	$$
	Then, if $\path$ is any discrete path and $\iota$ and admissible $\dist_\Mani$-marking of parameter $\delta_{n,1}=2\delta_{n}^{\dist_E}$, then $\iota$ is also and admissible $\dist_E$-marking of parameter
    $\delta_n^{\dist_E}$. Since $\dist_E(q_i,q_{i+1})^{\bar\alpha}\leq \dist_\Mani(q_i,q_{i+1})^{\bar\alpha}$ for any $\bar\alpha>0$, we obtain $C^{\dist_E}_n(x,y)\leq C^{\dist_\Mani}_{n,1}(x,y)$ for all $x,y\in\Mani$ and \eqref{eq:equivalence1} follows.
    \par
Let us now prove \eqref{eq:equivalence2}. Recall that $\delta_{n,2}=\delta_n^{\dist_E}$. 
Note that for any $q,q'\in\Mani$, $\dist_E(q,q')\geq\delta_n$ implies $\dist_\Mani(q,q')\geq\delta_n$. Thus
if $\iota$ is a $\dist_E$-admissible marking of a path $\path$ of parameter $\delta_n$, then it is also $\dist_\Mani$-admissible marking of the same parameter.
Furthermore, by Lemmas \ref{le:JnAsmp} and \ref{le:DeltanBound}, we may work on the event $\{J_n^{\dist_E}\leq \Delta_n/\theta\leq \delta_n\}$ for the $C_n^{\dist_E}$-geodesic path, say $\geod$. In this set, and using Proposition \ref{prop:euclideaVSgeodesic}, we have $\abs{\geodi_i - \geodi_{i+1}}\geq  \frac{\dist_\Mani(\geodi_i,\geodi_{i+1})}{1+\zeta_0(\delta_n)}$ for any
consecutive points $(\geodi_i,\geodi_{i+1})\subset\geod.$ As a consequence 
$C_n^{\dist_E}(x,y)\geq \frac{C_n^{\dist_\Mani}(x,y)}{(1+\zeta_0(\delta_n))^{\alpha+\frac{\maxv A}{\log n}}}$. Since
$(1+\zeta_0(\delta_n))^{\alpha+\frac{\maxv A}{\log n}}\toas{n\to\infty}1$, and the event
$\{J_n^{\dist_E}\leq \Delta_n/\theta\leq \delta_n\}$ has probability at least $1-e^{-\const_0n^{\const_1}}$ for some positive $\const_0,\const_1$, this proves \eqref{eq:equivalence2}.
The proof of the proposition is therefore concluded. 
\end{proof}

\section{Proofs of auxiliary results}\label{sec:appendix}
In this appendix we prove the technical lemmas stated in Section \ref{sec:prelim}. 

\subsection{Proofs of Lemmas  \ref{le:JnAsmp} and \ref{le:preliminaryUpperBoundSFD}}\label{sec:proofs_conddensityf_DeltanBound_JnAsmp_preliminaryUpperBoundSFD}

\begin{proof}[Proof of Lemma \ref{le:JnAsmp}]
Recall the function $\zeta$ of Proposition \ref{prop:euclideaVSgeodesic}.
Fix $0<\theta<1/2$ satisfying
$$
\forall n\geq2,\quad 0<2\left(\frac{1}{2}+\theta\right)^\alpha (1+\zeta(\theta))^{\alpha+\frac{\maxv A}{\log n}} <1.
$$
Note that such a $\theta$ can always be found, and be uniform on $n$, since $2^{-1/\alpha}>2^{-1}$ whenever $\alpha>1$, and also $\maxv A/\log n\toas{n\to\infty}0$ and, by Proposition \ref{prop:euclideaVSgeodesic},
 $\lim_{\theta\to0}\zeta(\theta)=0.$\par

Let $(\alpha_i)_{i\geq1}$ be as in \eqref{eq:defDFermatMicroGeneral}. For any $z,z'\in\Mani$ define the set
$$
\Theta_{i}(z,z')\coloneqq \{u\in \Mani\colon \dist_\Mani(z,u)^{\alpha_i} + \dist_\Mani(u,z')^{\alpha_i}  \leq \abs{z-z'}^{\alpha_i}\}.
$$ 
Note that we always have $\abs{z-z'}\leq \dist_\Mani(z,z').$ 
Let us assume that the pair $z,z'$ satisfies 
\begin{align} \label{eq:Jnc1}
\dist_\Mani(z,z')<\theta<1. 
\end{align}
If $z$ and $z'$ do not satisfy  \eqref{eq:Jnc1}, the following proof will still hold if we instead consider $\Theta_i'(z,z'')=\{u\in \Mani\colon \dist_\Mani(z,u)^{\alpha_i} + \dist_\Mani(u,z'')^{\alpha_i}  \leq \abs{z-z''}^{\alpha_i}\}$ where $z''$ is chosen along a $\dist_\Mani$-geodesic curve between $z$ and $z'$ to satisfy $\dist_\Mani(z,z')=\dist_\Mani(z,z'')+\dist_\Mani(z'',z')$, and also taking   $z''$ close enough to $z$ to ensure
$\dist_\Mani(z,z'')<\theta$. For simplicity in the notation let us directly assume \eqref{eq:Jnc1}.\par
Note that \eqref{eq:Jnc1} implies
$$
u \in \Theta_i(z,z') \implies \quad \dist_\Mani(z,u)<1 \quad \&\quad \dist_\Mani(u,z')<1. 
$$
Then, if  a sample point $q\in Q_n$ falls in the set $\Theta_i(z,z')$, we have
% for $j\in\{i,i+1\}$ and using that  $\alpha_{j}\geq \alpha_i$:
\begin{align*}
	\dist_\Mani(z,z')^{\alpha_i}\geq \abs{z-z'}^{\alpha_i}
	&\geq \dist_\Mani(z,q)^{\alpha_i} + \dist_\Mani(q,z')^{\alpha_{{i}}} 
	\\&\geq \abs{z-q}^{\alpha_i}+\abs{q-z'}^{\alpha_{{i}}}.
\end{align*} 
In particular,  if
$\Theta_i(z,z')\cap Q_n\neq\emptyset$ for every $i\geq 1$, then the points
$(z,z')$ cannot be contiguous points
in any $C_n$-geodesic  $\geod_n(x,y)$ for any $x,y\in\Mani$. The main argument of the proof consists in showing that indeed $\card\{\Theta_i(z,z')\cap Q_n\}\geq1$ for every $i\geq 1$ as soon as
$\dist_\Mani(z,z')\geq \Delta_{n}/\theta$.\par

Note that $z,z'$ being connected by a continuous $\dist_\Mani$-geodesic path,  there exists  $m_{z,z'}\in\Mani$ such that  
$\dist_\Mani(z,m_{z,z'})=\dist_\Mani(z,z')/2$ and $\dist_\Mani(z',m_{z,z'})=\dist_\Mani(z,z')/2$.
We now prove that 
\begin{equation}\label{eq:ballinTheta}
 \ball_\Mani(m_{z,z'},\theta\dist_\Mani(z,z')) %\subset \left\{\ball_{\R^D}(m_{z,z'},\theta(1+\zeta(\theta_0))\abs{z-z'})\cap \Mani\right\}
 \subset\Theta_i(z,z').
\end{equation}

Indeed, for any $y\in \ball(m_{z,z'},\theta\dist_\Mani(z,z'))$ we have, using triangle inequalities:
$$
\dist_\Mani(z,y)^{\alpha_i} + \dist_\Mani(y,z')^{\alpha_i}
\leq \dist_\Mani(z,z')^{\alpha_i}(1/2+\theta)^{\alpha_i} + \dist_\Mani(z,z')^{\alpha_{i}} (1/2+\theta)^{\alpha_{i}}.
$$
Since $\alpha_i\leq \alpha+\maxv A/\log n$, the choice of $\theta$ and Proposition \ref{prop:euclideaVSgeodesic} together with \eqref{eq:Jnc1}  yield
\begin{align*}\label{eq:ThetaBound}
&\dist_\Mani(z,z')^{\alpha_i}(1/2+\theta)^{\alpha_i} + \dist_\Mani(z,z')^{\alpha_{i}} (1/2+\theta)^{\alpha_{i}}
\nonumber\\
&\le \dist_\Mani(z,z')^{\alpha_i}2(1/2+\theta)^{\alpha} \leq \abs{z-z'}^{\alpha_i}2(1/2+\theta)^{\alpha}(1+\zeta(\theta))^{\alpha_i}
\\&<\abs{z-z'}^{\alpha_i}.
\end{align*}
This proves \eqref{eq:ballinTheta}.\par
By definition of $\Delta_n$ (in \eqref{eq:defSpacing}),
if $\dist_\Mani(z,z')>\Delta_n/\theta\ge \Delta_n$ then $\card\{Q_n\cap\ball(m_{z,z'},\theta\dist_\Mani(z,z'))\} \asgeq 1,$ so that \eqref{eq:ballinTheta} 
implies
$\card\{Q_n\cap\Theta_i(z,z')\} \asgeq 1$ for all $i\geq1$, and the proof is finished.
\end{proof}

%%%%%%%%%%%%%%%%%%%%%%%%%%%%%%%%%%%%%%%%%%%%%%%%%%%%%%%%%%%%%%%%%%%%%%%%%%%%%%%%%%%%

\begin{proof}[Proof of Lemma \ref{le:preliminaryUpperBoundSFD}]
 By Lemma \ref{le:JnAsmp}, the first inequality $C_n(x,y)\leq D_{n,\alpha}(x,y)$ is a simple consequence of the fact that for any $x\leq1$ and $a>0$ we have
    $x^{\alpha + a}\leq x^{\alpha}.$\par
    By Lemma \ref{le:DeltanBound} we may work on the set $\{\Delta_n/\theta\leq 1\}.$ Then,  \eqref{eq:uppberBoundCn} follows  from \eqref{eq:damelin_concentration} which yields, for sufficiently large $n$,
    \begin{equation}
%    \label{eq:damelin_concentration}
 \P\left(\sup_{\substack{\forall x,y\in\Mani\colon \\
                  \dist_\Mani(x,y)\geq \delta_n}}  {\frac{D_{n,\alpha}(x,y)}{n^{-\frac{\alpha-1}{d}} D_{\infty,\alpha}(x,y)} } \geq \mu(d,\alpha) +\eps \right)\leq e^{-\const_0 (n\delta_n^d)^\alpha},
    \end{equation}
    where $\const_0>0$ depends on $\eps.$
\end{proof}

%%%%%%%%%%%%%%%%%%%%%%%%%%%%%%%%%%%%%%%%%%%%%%%%%%%%%%%%%%%%%%%%%%%%%%%%%%%%%%%%%%%%

\subsection{Bounding path sizes via path costs.}
The following lemma is needed to prove Proposition \ref{prop:KnAsmp} in the next subsection.
It adapts the arguments in Lemma 11 \cite{HwangDamelinHero2016} into
a more general setting, allowing to prove the tighter upper bound for $K_n(x,y)$ given by Proposition \ref{prop:KnAsmp}.
Let us introduce the notation
\begin{equation}\label{eq:dA}
    \dist_\alpha(\path)=\sum_{i=1}^{\card\path-1}\dist(q_i,q_{i+1})^\alpha.
\end{equation}
\begin{lemma}\label{le:HDHLemma11}
Let $q,q'\in Q_n$ satisfy $\dist_\Mani(q,q')\geq \delta_n/2$. 
Let $\const_a>0$ be arbitrary and set
\begin{align*}
	a_n(q,q')&=\const_aD_{\infty,\alpha}(q,q'),\\
	c_n(q,q')&=a_n(q,q')n^{-\frac{\alpha-1}{d}}
\end{align*}
and
$$
 \ell_n(q,q')\coloneqq \lfloor \const_{K}a_n(q,q')n^{\frac{1}{d}}\rfloor
$$
where 
\begin{equation}\label{eq:constantNumberSteps}
	\const_{K}\coloneqq (2\const_d\maxv{f}\Gammaf(1+d/\alpha)e^{1+d/\alpha})^{\frac{\alpha+1}{d}}\frac{\alpha+2}{d}e^{\frac{\maxv A}{d}} + 1. 
\end{equation}
\par
Letting $$\alpha_n = \alpha + \frac{\maxv A}{\log(n)}$$  
we have, for large enough $n\geq n_0$, 
\begin{align}\label{eq:KnUpp_probBoundStatement}
&\P\left(\exists \path,q,q'\colon  q\qsim q',\left\{\path,q,q'\right\}\subset Q_n,\quad \dist_{\alpha_n}(\path)\leq c_n(q,q'), \quad
\card\path+2\ge\ell_n(q,q') \right)\\
&\leq e^4 n^2\exp\left\{-\frac{\const_K\const_a \minv f^{\frac{\alpha-1}{d}}}{2C}n^{\frac{\eta}{d}}
                    \right\}\nonumber
\end{align}
where we recall the constant $C>0$ from Assumption \ref{cond:DeltanBound} and notation $\dist_\alpha$ was defined in \eqref{eq:dA}.
\end{lemma}
\begin{remark}\label{rem:pathSizeBound1}
Looking at the proof of  Lemma \ref{le:HDHLemma11} it is clear that the same result holds if we replace $\maxv A$ by some $a>\maxv A$, changing the constant $\const_K$ accordingly and also possibly $n_0$.
\end{remark}

\begin{proof}[Proof of Lemma \ref{le:HDHLemma11}]
    Let $x,y\in\Mani$. In the following $\alpha>1$, 
    $\ell\geq1$, and $0<r\leq 1$ are arbitrary values which, when convenient, will be replaced
    by $\alpha=\alpha_n, \ell=\ell_n(x,y)-2,$ and $r=c_n(x,y)$. Define the set of discrete paths
    $$\{\path\colon x\qsim q, q\in Q_n, \dist_{\alpha}(\path)\leq r, \card\path \eq  \ell \}.$$
    That is the set of $Q_n$-paths starting at $x\in\Mani$ and ending in some $q\in Q_n$, whose weighted path length is less than
    $r$, and whose number of elements is exactly $\ell$. 
    Note that
    \begin{align*}
    &\left\{\exists \gamma,q\colon x\qsim q, q\in Q_n, \dist_{\alpha}(\path)\leq r, \card\path\geq \ell\right\} 
    \\&\subset\left\{\exists \gamma,q\colon x\qsim q, q\in Q_n, \dist_{\alpha}(\path)\leq r, \card\path\eq \ell\right\}.
    \end{align*}
    We will show that for large enough $n$, and $\forall\substack{x,y\in\Mani\\\dist_\Mani(x,y)\geq \delta_n/2}$, we have
    \begin{align}\label{eq:Knupperbound_prob1}
    &\P\left(
    \left\{\exists \gamma,q\colon x\qsim q, q\in Q_n, \dist_{\alpha_n}(\path)\leq c_n(x,y), \card\path\eq \ell_n(x,y)-2\right\}
    \right)\nonumber\\ 
    &\leq \exp\left\{-\ell_n(x,y)+3\right\}.
    \end{align}
    Note that since $\dist_\Mani(x,y)\geq \delta_n/2$ 
    we have 
    \begin{equation}\label{eq:anLowerBound}
    a_n(x,y)\geq \const_a \minv f^{\frac{1-\alpha}{d}}\delta_n/2.
    \end{equation}
    Then, using  Assumption \ref{cond:DeltanBound} we obtain, 
    \begin{equation}\label{eq:Knupperbound_barKnbound}
    \forall \substack{x,y\in\Mani\\ \dist_\Mani(x,y)\geq \delta_n/2}\colon \quad \ell_n(x,y) \geq \lfloor\const n^{\frac{\eta}{d}}\rfloor \uparrow \infty,     
    \end{equation}
    where $\const=\const_K\const_a \minv f^{\frac{\alpha-1}{d}}/2.$ Thus \eqref{eq:Knupperbound_prob1} is bounded by
    \begin{align*}
    &\P\left(
    \left\{\exists \gamma,q\colon x\qsim q, q\in Q_n, \dist_{\alpha_n}(\path)\leq c_n(x,y), \card\path\eq \ell_n(x,y)-2\right\}
    \right)\\ 
    &\leq e^4e^{- \const n^{\frac{\eta }{d}}}
    \end{align*}
    for large enough $n$.
    The stated result \eqref{eq:KnUpp_probBoundStatement} will then follow by summing the above over the possible choices of $x,y\in Q_n$.\par
    Let us then prove \eqref{eq:Knupperbound_prob1}. 
    To ease notation let us consider $\ell\in\N$. For any choice of
    $\alpha>1,\ell\geq 1, 0<r<1,$ define the random variable
    $$
    F_{\alpha,\ell,Q_n}(x,r)=\card \{\path\colon x\qsim q, q\in Q_n, \dist_{\alpha}(\path)\leq r, \card\path\eq \ell+1 \}.
    $$
    Let $\hat \ell_n(x,y)=\ell_n(x,y) -3$.
    By Markov's inequality eq. \eqref{eq:Knupperbound_prob1} will follow once we prove that for large enough $n$,
    \begin{equation}\label{eq:KnBoundNumberPaths}
    \forall\substack{x,y\in\Mani\\\dist_\Mani(x,y)\geq \delta_n/2}\colon\quad\E[F_{\alpha_n,\hat \ell_n(x,y),Q_n}(x,c_n(x,y))] \leq \exp\left\{-\hat \ell_n(x,y)\right\}.
    \end{equation}
    To prove the latter we use the recursive relation 
    \begin{align*}
    F_{\alpha,\ell,Q_n}(x,r) &= \sum_{q\in Q_n} \Ind{}\{\dist_\Mani(x,q)<r^{1/\alpha}\}F_{\alpha,\ell-1,Q_n\setminus \{q\}}(q, r - \dist_\Mani(x,q)^{\alpha});
    \end{align*}
   and $F_{\alpha,1,Q_n}(x,r)=\card\{Q_n \cap \ball(x,r^{1/\alpha}\}.$  %Using this bound, together with the Campbell-Mecke formula (Theorem 3.2 and Example 3.3 in \cite{BaddeleyBaranySchneider2006})
    Since, for every $q\in Q$, the set $Q\setminus \{q\}$ is i.i.d. we also obtain the following recursion in expectation
    \begin{equation}\label{eq:campMeck}
    \E[F_{\alpha, \ell,Q_n}(x,r)] = n\int_{\ball(x,r^{1/\alpha})}\mu(dq) \E[F_{\alpha,\ell-1,Q_{n-1}}(q, r -  \dist_\Mani(x,q)^{{\alpha}})];
    \end{equation}
    and also $\E[F_{\alpha,1,Q_n}(x,r)]= n\int_{\ball(x,r^{1/\alpha})}\mu(dq) 1$.\par 
    Let $r_0=r$. Set also $x_0=x$, and for a sequence of points $(x_1,\cdots,x_\ell)$
    set
    $r_k\equiv r_k(x_0,\cdots,x_k)\coloneqq r_0 - \sum_{i=1}^k \dist_\Mani(x_{i-1},x_{i})^{\alpha}$. Then we have, applying \eqref{eq:campMeck} recursively, 
    \begin{align*}
    &\E[F_{\alpha,\ell,Q_n}(x_0,r_0)]
    \\&= n(n-1)\cdots (n-\ell+1) \int_{\ball(x_0,r_0^{1/\alpha})}\mu(dx_1) \cdots \int_{\ball(x_{\ell-1},r_{\ell-1}^{1/\alpha_{n}})}\mu(dx_{\ell})\quad 1.\nonumber
    \end{align*}
    By \eqref{eq:density} we have, for the rightmost integral,
   
    $$
     \mu( \ball(x_{\ell-1},r_{\ell-1}^{1/\alpha}))\leq \const_d\maxv f r_{\ell-1}^{d/\alpha} (1+\lo(1)),
    $$
    where the error term $\lo(1)$ holds as $r_{\ell-1}\leq r_0\downarrow 0$ and uniformly on $x_{\ell-1}$. Plugging in the above we obtain
    \begin{align}\label{eq:campMeck2}
    &\E[F_{\alpha,\ell,Q_n}(x_0,r_0)]\\
    &\leq n^{\ell} \int_{\ball(x_0,r_0^{1/\alpha})}\mu(dx_1) \cdots \int_{\ball(x_{\ell-2},r_{\ell-2}^{1/\alpha})}\mu(dx_{\ell-1})\const_d\maxv{f} (1+\lo(1)) r_{\ell-1}^{d/\alpha}\nonumber 
    \end{align}
    for small enough $r_0$. 
    Let $\beta\geq0$. Recall that $r_{\ell-1}=r_{\ell-2}-\dist_\Mani(x_{\ell-1},x_{\ell-2})^{\alpha}<r_{\ell-2}$ so that 
    $\mu(r_{\ell-1}^\beta>s)=0$ whenever $s\geq r_{\ell-2}^{\beta}$. Thus, making the
    change of variable $u=s/r_{\ell-2}^{\beta}$ below, we compute
    \begin{align}\label{eq:campMeck3}
    %&\int_{\ball(x_{n-2},g_{n-1}^{-1}(r_{n-2}))}\mu(dx_{n-1})\left(g_{n}^{-1}(r_{n-1})\right)^{d}\nonumber \\
    &\int_{\ball(x_{\ell-2},r_{\ell-2}^{1/\alpha})}\mu(dx_{\ell-1})r_{\ell-1}^{\beta}
    =\int_0^{r_{\ell-2}^{\beta}} \mu\left( x_{\ell-1}\colon r_{\ell-1}^{\beta} > s \right) ds\nonumber\\
    &=r_{\ell-2}^{\beta} \int_0^{1} \mu\left( x_{\ell-1}\colon  r_{\ell-1}^{\beta} > u r_{\ell-2}^{\beta} \right) du\nonumber\\
    &=r_{\ell-2}^{\beta} \int_0^{1} \mu\left( x_{\ell-1}\colon r_{\ell-2}-\dist_\Mani(x_{\ell-1},x_{\ell-2})^{\alpha}> u^{1/\beta} r_{\ell-2} \right) du\nonumber\\
    &=r_{\ell-2}^{\beta} \int_0^{1} \mu\left( x_{\ell-1}\colon \dist_\Mani(x_{\ell-1},x_{\ell-2}) < \left(r_{\ell-2}- u^{1/\beta} r_{\ell-2}\right)^{1/\alpha} \right) du.
    \end{align}
    Using \eqref{eq:density}, the r.h.s. of \eqref{eq:campMeck3} is bounded by
    \begin{align*}
    &\const_d\maxv{f}(1+\lo(1))r_{\ell-2}^{\beta} \int_0^{1}  \left(r_{\ell-2}- u^{1/\beta} r_{\ell-2}\right)^{d/\alpha} du\\
    &=\const_d\maxv{f}(1+\lo(1))r_{\ell-2}^{\beta}r_{\ell-2}^{d/\alpha} \int_0^{1}  \left(1- u^{1/\beta} \right)^{d/\alpha} du
    \end{align*}
    where the error term $\lo(1)$ is again uniform on $x_{\ell-2}$ as $r_0\downarrow 0$. Making the change of variable $s=u^{1/\beta}$ in the last integral above we obtain
    \begin{align*}
     \int_0^{1}  \left(1- u^{1/\beta} \right)^{d/\alpha} du = \beta \int_0^1 (1-s)^{d/\alpha}s^{\beta -1}du
    =\frac{\Gammaf(1+d/\alpha)\Gammaf(1+\beta)}{\Gammaf(1+d/\alpha+\beta)}.
    \end{align*}
    Plugging in \eqref{eq:campMeck3} with $\beta\equiv \beta_1=d/\alpha$ we obtain
    \begin{align*}
    &\int_{\ball(x_{\ell-2},r_{\ell-2}^{1/\alpha})}\mu(dx_{\ell-1})r_{\ell-1}^{d/\alpha} \\&\leq (1+\lo(1))\const_d\maxv{f}r_{\ell-2}^{d/\alpha+d/\alpha} \frac{\Gammaf(1+d/\alpha)\Gammaf(1+d/\alpha)}{\Gammaf(1+d/\alpha+d/\alpha)}.
    \end{align*}
    Letting $\beta_k=dk/\alpha$, using the above computation recursively over $\beta_2,\cdots,\beta_\ell$ in 
     \eqref{eq:campMeck2}, and canceling terms, we arrive at
    \begin{align}\label{eq:Knupperbound_moment1}
        \E[F_{\alpha,\ell,Q_n}(x_0,r_0)] \leq (1+\lo(1))^\ell  \const_d^{\ell}\maxv{f}^\ell n^\ell r_0^{d\ell/\alpha} 
         \frac{ \Gammaf(1+d/\alpha)^\ell }{\Gammaf(1+\frac{d\ell}{\alpha})}
    \end{align}
    where the error term $\lo(1)$ holds as $r_0\downarrow0$ uniformly on $x_0\in\Mani$. 

    Set $\ell=\hat\ell_n\equiv \hat\ell_n(x,y)$, $r_0=c_n\equiv c_n(x,y)$, and $\alpha=\alpha_n$. 
    Note that since the sequence $(\alpha_n)_{n\geq1}$ is bounded from above, and also using \eqref{eq:Knupperbound_barKnbound}, by Stirling's approximation we have, as 
    $n\to\infty$ and uniformly on $\{x,y\in\Mani\colon \dist_\Mani(x,y)\geq\delta_n/2\}$,
    $$
     \Gammaf(1+\frac{d\hat\ell_n(x,y)}{\alpha_n}) = e^{-\frac{d\hat\ell_n(x,y)}{\alpha_n}}(\frac{d\hat\ell_n(x,y)}{\alpha_n})^{\frac{d\hat\ell_n(x,y)}{\alpha_n}} \bO(\sqrt{\hat\ell_n(x,y)}).
    $$
    Plugging in \eqref{eq:Knupperbound_moment1} we obtain
    \begin{align}\label{eq:Knupperbound_moment2}
    \forall \substack{x,y\in\Mani\\\dist_\Mani(x,y)\geq \delta_n/2}\colon\quad &\E[F_{\alpha_n,\hat\ell_n(x,y),Q_n}(x,c_n(x,y))]  
    \nonumber
    \\&\leq 
    \left(\frac{(1+\lo(1))  \const_d\maxv{f}\Gammaf(1+d/\alpha_n)e^{d/\alpha_n}n c_n(x,y)^{d/\alpha_n}}
    {(\frac{d\hat\ell_n(x,y)}{\alpha_n})^{d/\alpha_n}}\right)^{\hat\ell_n(x,y)}
    \nonumber\\&\times \bO(1/\sqrt{\hat\ell_n(x,y)}).
    \end{align}
    The error term $\lo(1)$ holds as $c_n(x,y)\to 0$, and since $c_n(x,y)\to 0$ as $n\to\infty$ uniformly on $x,y\in\Mani$,
    the error $\lo(1)$ holds as $n\to\infty$ uniformly on $x,y\in\Mani$.\par
    The quantity in the r.h.s. of \eqref{eq:Knupperbound_moment2}  decays exponentially in $\hat\ell_n(x,y)$ as soon as
    % $$
    % \liminf_{n\to\infty} \frac{(1+\lo(1))\maxv{f}\Gammaf(1+d/\alpha_n)e^{\phi(n)}}{\phi(n)^{\phi(n)}} \frac{n r_n^{\phi(n)}}{\ell_n^{\phi(n)}} < 1, 
    % $$
    % or
    $$
    \frac{(1+\lo(1))  \const_d\maxv{f}\Gammaf(1+d/\alpha_n)e^{d/\alpha_n}n c_n(x,y)^{d/\alpha_n}}
    {(\frac{d\hat\ell_n(x,y)}{\alpha_n })^{d/\alpha_n }} \leq \frac{1}{e}.
    $$
    We may choose $n\geq n_0$ large enough so that the error term $\lo(1)$ is less than $1$ (uniformly on $x,y\in\Mani$). Then the above inequality is satisfied whenever
    $$
   \left(\frac{d\hat\ell_n(x,y)}{\alpha_n +1}\right)^{d/\alpha_n }
   \geq 2\const_d\maxv{f}\Gammaf(1+d/\alpha_n)e^{1+d/\alpha_n}n c_n(x,y)^{d/\alpha_n};
    $$
    i.e. whenever
    $$
    \hat\ell_n(x,y) \geq (2\const_d\maxv{f}\Gammaf(1+d/\alpha_n)e^{1+d/\alpha_n})^{\frac{\alpha_n}{d}}\frac{\alpha_n+1}{d}n^{\frac{\alpha_n}{d}} c_n(x,y).
    $$  
    Taking large enough $n\geq n_0$
    so that $\maxv A/\log(n) <1$, the latter holds in particular whenever
    \begin{align}\label{eq:ellnLowerBound1}
    &\ell_n(x,y) 
    \nonumber\\&\geq \left\lfloor \left((2\const_d\maxv{f}\Gammaf(1+d/\alpha)e^{1+d/\alpha})^{\frac{\alpha+1}{d}}\frac{\alpha+2}{d}e^{\frac{\maxv A}{d}} + \frac{2n^{-\frac{1}{d}}}{a_n(x,y)}  \right) a_n(x,y)n^{\frac{1}{d}}\right\rfloor.
    \end{align}
    By \eqref{eq:anLowerBound} and Assumption \ref{cond:DeltanBound} we have $a_n(x,y)n^{\frac{1}{d}}\geq  \frac{\const_a \minv f^{\frac{1-\alpha}{d}}}{2C}n^{\frac{\eta}{d}}\toas{n\to\infty}\infty$ uniformly on $\{x,y\in\Mani\colon \dist_\Mani(x,y)\geq\delta_n/2\}$. Thus we may choose
    $n$ large enough to ensure that \eqref{eq:ellnLowerBound1} is satisfied whenever
    $$
    \ell_n(x,y) \geq \left\lfloor \left((2\const_d\maxv{f}\Gammaf(1+d/\alpha)e^{1+d/\alpha})^{\frac{\alpha+1}{d}}\frac{\alpha+2}{d}e^{\frac{\maxv A}{d}} + 1  \right) a_n(x,y)n^{\frac{1}{d}}\right\rfloor.
    $$
    Equation \eqref{eq:KnBoundNumberPaths} then follows by taking $\const_K$ as in \eqref{eq:constantNumberSteps}.
\end{proof}

%%%%%%%%%%%%%%%%%%%%%%%%%%%%%%%%%%%%%%%%%%%%%%%%%%%%%%%%%%%%%%%%%%%%%%%%%%%%%%%%%%%%
\subsection{Proof of Proposition \ref{prop:KnAsmp}}\label{sec:proof_prop:KnAsmp}

\begin{proof}[Proof of Proposition \ref{prop:KnAsmp}]
    The first statement follows from the simple bound
    $$
    K_n J_n \asgeq \sum_{i=1}^{K_n-1} \dist_\Mani(\geodi_{i+1},\geodi_i) \asgeq \dist_\Mani(x,y)
    $$
    so that, by Lemma \ref{le:JnAsmp},
    $$
    K_n \asgeq \dist_\Mani(x,y)/J_n \asgeq \dist_\Mani(x,y)\theta/\Delta_n.
    $$\par
On the other hand, the high-probability upper-bound for $K_n$ follows from Lemmas \ref{le:preliminaryUpperBoundSFD} and \ref{le:HDHLemma11}. Indeed, for any $x,y\in\Mani$ define
$$\path(x,y)\equiv (\pathi_2,\cdots,\pathi_{K_n-1})\equiv (\pathi_2(x,y),\cdots,\pathi_{K_n-1}(x,y))$$
to be the optimal path $\geod(x,y)$ except that the start and end points ($x$ and $y$) are removed.
By Lemma \ref{le:DeltanBound} we may work on the high-probability set $\{\Delta_n/\theta\leq 1\}$. In this set,
using Lemma \ref{le:JnAsmp}, we have that any pair of consecutive elements in $\geod(x,y)$ are at a distance smaller than $1$. Then,
removing the cost of the first and last jump in $\geod(x,y)$, and also bounding $A(s)/\log(n) \leq \maxv A/\log(n)$,
we obtain
\begin{align}\label{eq:KnLower_pathCostUpper}
 \dist_{\alpha + \frac{\maxv A}{\log(n)}}(\path(x,y))
 %\dist_{A}(\geod(x,y))=
 %\red{\sum_{i=2}^{K_n(x,y)-2} \dist(\pathi_{i+1},\pathi_i)^{\alpha + \frac{\maxv A}{\log(n)}}}
 \leq
C_n(x,y).
\end{align}\par
 Moreover,
by  Lemma \ref{le:JnAsmp} and
on the high-probability set $\{\Delta_n/\theta\leq \delta_n/4 \}$ we obtain that, whenever
$\dist_\Mani(x,y)\geq\delta_n$, we also have for the corresponding path $\path(x,y)$, that 
\begin{equation}\label{eq:Knlower_distqqprime}
\forall \substack{x,y\in\Mani\\ \dist_\Mani(x,y)\geq\delta_n}\colon\quad \dist_\Mani(\pathi_2(x,y),\pathi_{K_n-1}(x,y))\geq \delta_n/2.
\end{equation}\par
On the other hand,
\begin{align*}
\abs{D_{\infty,\alpha}(\pathi_2,\pathi_{K_n-1})- D_{\infty,\alpha}(x,y)}\leq 2\minv f^{-\frac{\alpha-1}{d}}\Delta_n/\theta.
\end{align*}
Then, using \eqref{eq:Knlower_distqqprime} and still on the set $\{\Delta_n/\theta\leq \delta_n/4 \}$ we have that $\forall \substack{x,y\in\Mani\\ \dist_\Mani(x,y)\geq\delta_n}$:
\begin{align*}
\max\left\{\abs{\frac{D_{\infty,\alpha}(\pathi_2,\pathi_{K_n-1})}{D_{\infty,\alpha}(x,y)}-1},
\abs{\frac{D_{\infty,\alpha}(x,y)}{D_{\infty,\alpha}(\pathi_2,\pathi_{K_n-1})}-1}\right\}
\leq  \frac{2\minv f^{-\frac{\alpha-1}{d}}\Delta_n/\theta}{\maxv f^{-\frac{\alpha-1}{d}}\delta_n/2}.
\end{align*}
Thus, on the set
$\left\{\frac{2\minv f^{-\frac{\alpha-1}{d}}\Delta_n/\theta}{\maxv f^{-\frac{\alpha-1}{d}}\delta_n/2}\leq 1\right\}\subset \{\Delta_n/\theta\leq \delta_n/4 \}$ we obtain
\begin{align}\label{eq:KnUpper_DxyLowerBound}
\forall \substack{x,y\in\Mani\\ \dist_\Mani(x,y)\geq\delta_n}\colon\quad \max\left\{\frac{D_{\infty,\alpha}(\pathi_2,\pathi_{K_n-1})}{D_{\infty,\alpha}(x,y)},
\frac{D_{\infty,\alpha}(x,y)}{D_{\infty,\alpha}(\pathi_2,\pathi_{K_n-1})}\right\}
\leq  2.
\end{align}
From \eqref{eq:KnUpper_DxyLowerBound}  we conclude, on the one hand,
\begin{equation*}
\forall \substack{x,y\in\Mani\\ \dist_\Mani(x,y)\geq\delta_n}\colon\quad D_{\infty,\alpha}(\pathi_2,\pathi_{K_n-1})\geq D_{\infty,\alpha}(x,y)/2.
\end{equation*}
Take now $\const_a=4\mu(d,\alpha)$ and recall the corresponding functions $a_n(q,q')=\const_a D_{\infty,\alpha}(q,q')$, $c_n(q,q')=a_n(q,q')n^{-\frac{\alpha-1}{d}}$, and $\ell_n(q,q')$ from Lemma \ref{le:HDHLemma11}.
In particular note that
\begin{equation}\label{eq:Knupper_cnLowerBound1}
%\forall \substack{x,y\in\Mani\\ \dist_\Mani(x,y)\geq\delta_n}\colon\quad
c_n(\pathi_2,\pathi_{K_n-1})= \const_a D_{\infty,\alpha}(\pathi_2,\pathi_{K_n-1})n^{-\frac{\alpha-1}{d}}\geq 2\mu(d,\alpha)D_{\infty,\alpha}(x,y)n^{-\frac{\alpha-1}{d}}.
\end{equation}
Combining \eqref{eq:KnLower_pathCostUpper} and \eqref{eq:Knupper_cnLowerBound1} together with
Lemma \ref{le:preliminaryUpperBoundSFD} we obtain, for large enough $n$,
\begin{align}\label{eq:KnUpperbound1_Cnupperbound}
&\P\left(\sup_{\substack{x,y\in\Mani\\\dist_\Mani(x,y)\geq \delta_n}} \frac{\dist_{\alpha + \frac{\maxv A}{\log(n)}}(\path(x,y))}{c_n(q_2(x,y),q_{K_n-1}(x,y))}\geq  1 \bigcap\left\{\frac{2\minv f^{-\frac{\alpha-1}{d}}\Delta_n/\theta}{\maxv f^{-\frac{\alpha-1}{d}}\delta_n/2}\leq 1\right\}  \right)\nonumber\\
&\leq \P\left(\sup_{\substack{x,y\in\Mani\\\dist_\Mani(x,y)\geq \delta_n}} \frac{C_n(x,y)}{\mu(d,\alpha)n^{-\frac{\alpha-1}{d}}D_{\infty,\alpha}(x,y)}\geq  2\right)
\leq e^{-\const_0' (n\delta_n^d)^{\const_1'}},
\end{align}
for some $\const_0',\const_1'>0$.  
On the other hand, from \eqref{eq:KnUpper_DxyLowerBound} we also have
\begin{align*}
\forall \substack{x,y\in\Mani\\ \dist_\Mani(x,y)\geq\delta_n}\colon\quad\upb K_n(x,y)&=\lfloor\const_{K} 2\mu(d,\alpha)D_{\infty,\alpha}(x,y)n^{\frac{1}{d}}\rfloor\\
&\geq \lfloor\const_{K} 4\mu(d,\alpha)D_{\infty,\alpha}(q_2(x,y),q_{K_n-1}(x,y))n^{\frac{1}{d}}\rfloor\\
&= \ell_n(\pathi_2(x,y),\pathi_{K_n-1}(x,y)).
\end{align*}
In other words, the latter observations
together with Lemma \ref{le:DeltanBound} imply the existence of constants $\const_0'>0$ and $\const_1'>0$
such that, 
$
\P(\mathcal A_n^c) \leq e^{-\const_0'n^{\const_1'}}
$
for the set
\begin{align*}
\mathcal A_n &= \left\{\upb K_n(x,y) \geq \ell_n(\pathi_2(x,y),\pathi_{K_n-1}(x,y))  \right\}
    \\&\bigcap \left\{ \sup_{\substack{x,y\in\Mani\\\dist_\Mani(x,y)\geq \delta_n}} \frac{\dist_{\alpha + \frac{\maxv A}{\log(n)}}(\path(x,y))}{ c_n(\pathi_1(x,y),\pathi_{K_n-2}(x,y))}\leq 1 \right\}\\
    \\&\bigcap \left\{\inf_{\substack{x,y\in\Mani\\\dist_\Mani(x,y)\geq \delta_n}}\dist_\Mani(\pathi_2(x,y),\pathi_{K_n-1}(x,y))\geq \delta_n/2  \right\}.
\end{align*}
 The above yields
\begin{align*}
&\P\left(
\sup_{\substack{x,y\in\Mani\\\dist_\Mani(x,y)\geq \delta_n}} \frac{K_n(x,y)}{\upb K_n(x,y)} > 1
\right)
\\&\leq e^{-\const_0'n^{\const_1'}}
+\P\left(\sup_{\substack{x,y\in\Mani\\\dist_\Mani(x,y)\geq \delta_n}} \frac{\card{\path(x,y)} + 2}{\ell_n(\pathi_2(x,y),\pathi_{K_n-1}(x,y))}\ge  1 \bigcap 
\mathcal A_n\right).
%&\leq e^{-\const_0' (n\delta_n^d)^{\const_1'}}+ 
\end{align*}
By Lemma \ref{le:HDHLemma11}, the probability in the r.h.s. is also exponentially bounded.\par
  Finally, repeating the same argument but replacing $\maxv A$ by some $a>\maxv A$, we still obtain \eqref{eq:KnBoundNumberPaths} (see Remark \ref{rem:pathSizeBound1}). Thus the result holds uniformly on the choice of $A$ as long as $\maxv A\leq a$. 
\end{proof}

\section{Statement of AI use}

The research direction, mathematical content, and overall proof strategy of this paper are our own. 
However, AI assistance (Claude) played a limited but genuinely useful mathematical role, beyond pure copy-editing. In particular, while we were working out the micro-to-meso step, repeated exchanges with the AI helped us arrive at a sharper definition of the microscopic model, specifically a cleaner formulation of the discretised path-cost. The AI did not produce the proofs, but it helped us recognise that an earlier, more natural-looking definition was an obstruction, and pushed the formulation in the direction that ultimately went through. \par
Also, the ideas and the structure of the proofs in Section \ref{subsec:locpropcinfty} were our own. However, a first draft (wording) of the standard arguments used therein was generated with AI assistance, and then revised and polished by the authors.\par
Beyond this, the AI was used for standard editorial tasks, proofreading of the mathematical proofs, as well as general improvement of the writing style and grammar of the manuscript.

\bibliographystyle{jtbnew}
\bibliography{Math}
%\printbibliography
\end{document}